# Competitive Robust Dynamic Pricing in Continuous Time with Fixed Inventories

Terry L. Friesz[a,*], Changhyun Kwon[a], Tae Il Kim[a], Lifan Fan[b], Tao Yao[a]

[a] Department of Industrial and Manufacturing Engineering, The Pennsylvania State University, University Park, PA 16802, USA

[b] Business School, Guangzhou University, Guangzhou, 510006, China

**Abstract:** The problem of robust dynamic pricing of an abstract commodity, whose inventory is specified at an initial time but never subsequently replenished, originally studied by Perakis and Sood (2006) in discrete time, is considered from the perspective of continuous time. We use a multiplicative demand function to model the uncertain demand, and develop a robust counterpart to replace the uncertain demand constraint. The sellers' robust best response problem yields a generalized Nash equilibrium problem, which can be formulated as an equivalent, continuous-time quasi-variational inequality. We demonstrate that, for appropriate regularity conditions, a generalized robust Nash equilibrium exists. We show that the quasi-variational inequality may be replaced by an equivalent variational inequality, and use a fixed-point algorithm to solve the variational inequality. We also demonstrate how explicit time lags associated with price updating in real-world decision environments, as well as specific pricing decision rules, may be introduced to create a dual time scale formulation and the associated solutions computed. We illustrate, via numerical examples, how robust pricing based on our DPFI formulation offers generally superior and never inferior worst case performance compared to nominal pricing.

* Corresponding author. Tel.: +1 814 863 2445; Fax: +1 814 863 4745

Email address: tfriesz@psu.edu (T.L. Friesz), chkwon@buffalo.edu (C. Kwon), fanlifan320@gmail.com (L. Fan), taoyao@psu.edu (T.Yao)

# 1 Introduction

## 1.1 Background

We explore the dynamic pricing of an abstract commodity in an oligopolistic market. Each seller's inventory is specified at an initial time but never subsequently replenished. Every seller competes with each other by setting its prices. Perakis and Sood (2006) refer to such problems as dynamic pricing with fixed inventory (DPFI) and study them with demand uncertainty over a finite discrete time horizon. They employ the ideas from robust optimization and quasi-variational inequality to address this problem.

There is a specific and compelling motivation for the study of DPFI problems: they form a theoretical foundation for a broad class of revenue management problems. More specifically, they are an abstraction of service systems whose output for a given planning horizon cannot be held in inventory. As such, they capture the essence of actual airline, ferry, rail, transit, and other transportation pricing and service timing decisions.

Unlike Perakis and Sood (2006), we study a generic DPFI using a continuous-time perspective, which includes computation in continuous time, and we model demand uncertainty in a quite different way. As we shall show, our approach to formulating a generic DPFI problem allows us to easily establish sufficiency of the quasivariational inequality that is a necessary condition for each firm's best response problem, in contrast to the much less direct demonstration of sufficiency in Perakis and Sood (2006).

Moreover, the use of continuous time allows us to easily and naturally extend the analysis of fixed inventory problems to explicitly consider dual time scales; thereby, we are able to introduce explicit pricing rules that recognize the intrinsic time lag that characterizes price adjustments in actual decision environments. Such features are not addressed by Perakis and Sood (2006).

Extension of some parts of the Perakis and Sood (2006) model to infinite dimensional vector spaces is trivial; extension of other parts of their discrete time model to continuous time are not straight forward and require carefully crafted arguments. This paper provides the detailed analysis, including re-casting the problem as a robust infinite dimensional quasi-variational inequality, and showing a robust generalized Nash equilibrium exists. We present a continuous-time fixed point algorithm whose subproblems are linear quadratic optimal control problems to obtain the robust generalized Nash equilibrium. Several fixed inventory numerical examples are presented.

Since dynamic optimization usually occurs in continuous time, it makes our continuous-time model more realistic. And our continuous-time perspective makes use of the notion of an infinite dimensional quasi-variational inequality. To the best of our knowledge, the theory of continuous-time quasi-variational inequality that we study in this paper has not previously been applied to revenue management and production planning.

Many service provision problems may be conceptualized as DPFI-problems. The sale of seats on a plane making a single outbound flight, as might arise for chartered aircraft, is an example. As such DPFI-problems represent a kind of proto-service environment that allows the undistracted exploration of issues intrinsic to dynamic pricing. In the spirit of Perakis and Sood (2006), we make the following assumptions:

1. Perfect information. We assume that perfect information obtains at the outset about the structure of demand, the impact of price changes on demand, and the initial inventory.
2. Consumer choice and demand. We assume that demand for the output of each seller is a function only of current prices, and prices are the only factor that distinguishes products from one another.
3. Product. We assume there is a single product and that inventory must be zero at the terminal time. This assumption is consistent with the view that inventory is saleable for all $t \in [t_0, t_f)$ and is worthless at $t_f$.
4. Objectives. We further assume that sellers maximize the present value of their respective revenues by setting prices and do not employ any other type of strategies. Accordingly a generalized Nash equilibrium will describe the market of interest.

**1.2 Literature Review**

Dynamic pricing has been extensively researched in the literature of revenue management. Gallego and van Ryzin (1994) consider a monopolistic dynamic pricing problem with fixed inventory over a finite continuous-time horizon. For a family of exponential demand functions, the optimal pricing policy can be derived in closed form. For general demand function, they find an upper bound on the expected revenue and obtain the form of near-optimal policies. Feng and Xiao (2000) also address a continuous-time monopolistic dynamic pricing problem. They assume reversible price changes are allowed and demand follows a Poisson process. The problem is formulated as an intensity control model, and the optimal prices are obtained in closed form. Levin et al. (2010) study a dynamic pricing problem for a monopolistic firm in the presence of strategic consumers who know that pricing is dynamic and may adjust their purchases accordingly. They formulate the problem as a stochastic dynamic game, and demonstrate the existence of a unique subgame-perfect equilibrium pricing policy. For surveys on monopolistic dynamic pricing problem, see Weatherford and Bodily (1992), McGill and van Ryzin (1999), Elmaghraby and Keskinocak (2003), Bitran and Caldentey (2003), and Talluri and van Ryzin (2004).

Recently, there is a growing interest to consider competition in dynamic pricing problem. Gallego et al. (2006) study a sequential game and a repeated simultaneous game in a duopoly market. Both firms with fixed inventory sell products in both a forward market and a spot market. Granot et al. (2010) address a multi-period, no-replenishment, dynamic pricing problem under duopoly competition. They assume that every consumer visits only one retailer in any period. If the price is lower than the consumer's valuation, he will purchase the product. Otherwise he will visit the other retailer in the next period.

Kachani et al. (2004) and Kwon et al. (2009) consider dynamic pricing with fixed inventory in an oligopolistic setting integrating demand learning. Kachani et al. (2004) assume demand is a linear function and formulate the problem as a Mathematical Program with Equilibrium Constraints (MPEC). Kwon et al. (2009) model demand dynamics as a differential equation based on evolutionary game theory and learn the demand parameters based on Kalman filtering. Lin and Sibdari (2009) introduce a game-theoretic model to describe dynamic price competition between firms. They prove the existence of Nash equilibrium. Levin et al. (2009) address a dynamic pricing problem for oligopolistic firms facing strategic consumers. They prove that a unique subgame-perfect equilibrium exists and discuss the impact of strategic consumer behavior on revenue. Gallego and Hu (2009) present a choice-based dynamic pricing model for oligopolistic firms selling persihable inventories over a continuous finite horizon. They show the existence and uniqueness of the open-loop Nash equilibrium.

Demand uncertainty is another important aspect in our paper. Zabel (1970) considers two approaches to model demand uncertainty: a multiplicative demand model ($d(t) = \xi(t)u(p(t))$) and an additive demand model ($d(t) = u(p(t)) + \xi(t)$), where $d(t)$ is the demand at time $t$, $p(t)$ is the price at time $t$, $u()$ is a linear decreasing function of price, and $\xi(t)$ is the uncertainty factor at time $t$. Chen and Simchi-Levi (2004) also introduce two models for demand function: an additive demand function ($D_t(p_t, \beta_t) = D_t(p_t) + \beta_t$) and a nonadditive demand function ($D_t(p_t, \varepsilon_t) = \alpha_t D_t(p_t) + \beta_t$), where $\varepsilon_t = (\alpha_t, \beta_t)$ and $\alpha_t, \beta_t$ are two random variables with $E\{\alpha_t\} = 1$ and $E\{\beta_t\} = 0$. Perakis and Sood (2006) model uncertain demand as a function of price and an uncertainty factor which is a vector of demand parameters that can take any value from a given closed and convex uncertainty set. In this paper, we use the multiplicative demand function to model demand uncertainty.

There are several ways to deal with demand uncertainty in optimization problems. Robust optimization is one of them, which we will use in this paper. The earliest research on robust optimization dates back to Soyster (1973), who considers a linear optimization problem and assumes all uncertain parameters take the worst-case values from a convex set. Ben-Tal and Nemirovski (1998) consider a robust convex optimization problem with an ellipsoidal uncertainty set and show that the robust convex progam of some convex optimization problems, such as linear programming, quadratically constrained programming, semidefinite programming, can be exactly or approximately solved by polynomial time interior point algorithm. Ben-Tal and Nemirovski (1999) address a linear programming with uncertain data, and replace the uncertain linear programming by a robust counterpart. In particular, for an linear programming with ellipsoidal uncertainty set, the corresponding robust counterpart can be solved in polynomial time. Perakis and Sood (2006) consider a competitive dynamic pricing problem based on robust optimization approach over a finite discrete time horizon. Bertsimas and Thiele (2006) employ the ideas of robust optimization to address the optimal inventory control problem in a supply chain. Aghassi and Bertsimas (2006) apply robust optimization to a distribution-free model of games with incomplete information, and provide a robust-optimization equilibrium. Leung et al. (2007) propose a robust optimization model to solve a multi-site production planning problem with uncertain data.

### 1.3 Contributions and Organization of the paper

The contributions of this paper are the following:

1. We extend the discrete-time DPFI model of Perakis and Sood (2006) to continuous-time DPFI model. To the best of our knowledge, this is the first paper to use robust optimization and quasi-variational inequality to address the oligopolistic dynamic pricing problem with demand uncertainty over a finite continuous-time horizon.
2. We model uncertain demand as a multiplicative demand function, and restate the uncertain demand constraint by a robust counterpart, which leads the robust optimization problem become a deterministic optimization problem.
3. The sellers' best response problem yields a generalized Nash equilibrium problem (GNEP) which can be represented as a quasi-variational inequality (QVI). We prove the equivalence between GNEP and QVI, and the existence of a generalized robust Nash equilibrium.
4. Since efficient convergent algorithm for solving QVI is not available, we construct a variational inequality (VI) which can be efficiently solved by a fixed point algorithm, and show that any solution of the VI is a solution of the QVI.
5. Numerical examples show that a seller using our robust pricing policy can always obtain a profit with standard deviation 0, and can improve much better the worst-case performance compared to nominal pricing policy. For some distributions of the uncertainty factor, robust pricing policy can generate higher average profit than nominal pricing policy. We also find that upon a certain level of robust magnitude, reducing standard deviation would sacrifice a slight of average profit.

This article is organized as follows: Section 2 states the notation and the regularity conditions that will be used throughout this paper. Section 3 models the best response problem for each seller. In section 4, we articulate a robust counterpart of the best response problem. Section 5 formulates the generalized robust Nash equilibrium as a quasi-variational inequality. In section 6, we build a VI and show that any solution of the VI is a solution of the QVI. A fixed point algorithm is introduced to solve the VI. Section 7 provides some numerical examples. Section 8 summarizes our conclusions.

## 2 Notation and Regularity Conditions

We will use $\pi_s(t)$ to denote the price charged by seller $s \in S$ at time $t \in [t_0, t_f]$, where $S$ denotes the set of all sellers and $\pi_s \in L^2[t_0, t_f]$, $t_0 \in R^1_+$, $t_f \in R^1_{++}$, and $t_f > t_0$, where $L^2[t_0, t_f]$ is the space of square-integrable functions, a Hilbert space. For any subscript $s \in S$ and any given instant of time $t_1 \in [t_0, t_f]$, we stress that $\pi_s(t_1) \in R_+^1$ is a scalar. Let $\pi_{-s} = (\pi_g : g \in S \setminus s)$ be a column vector of non-own prices. We also use the notation $\pi = (\pi_s : s \in S) \in (L^2[t_0, t_f])^{|S|}$ to represent the column vector of prices that are the decision variables of the model to be constructed. We use the notation $h_s(\pi(t); \xi_s) : (L^2[t_0, t_f])^{|S|} \times R_+^1 \to H^1[t_0, t_f]$ to represent the observed demand for the output of each seller $s \in S$ corresponding to a specific realization of the vector of parameters $\xi_s \in U_s \subseteq (L^2[t_0, t_f])^{|m|}$, where $U_s$ is uncertainty set. $H^1[t_0, t_f]$ is the Sobolev space of functions that belong to $L^2[t_0, t_f]$ and have derivatives also belonging to $L^2[t_0, t_f]$; furthermore $H^1[t_0, t_f]$ is a Hilbert space. We let $D_s(t)$ represent the realized demand served by seller $s \in S$, and define column vectors

$$D = (D_s : s \in S) \in (L^2[t_0, t_f])^{|S|}$$

$$h = (h_s : s \in S) \in (H^1[t_0, t_f])^{|S|}$$

Realized demands will also be decision variables. Since realized demand must be less than or equal to observed demand, for each seller $s \in S$, we impose the constraint

$$D_s(t) \leq h_s(\pi(t); \xi_s) \quad \forall s \in S, t \in [t_0, t_f] \quad (1)$$

or equivalently

$$D(t) \leq h(\pi(t)) \quad \forall t \in [t_0, t_f] \quad (2)$$

Throughout this paper, unless otherwise noted, we invoke the following regularity conditions:

A1. Prices are bounded from above and below according to

$$\pi_s \geq \pi_{\min} \in R_+^1$$

$$\pi_s \leq \pi_{\max} \in R_{++}^1$$

for all $s \in S$, where $\sup_{\xi_s, \pi_{-s}} h_s(\pi_{\max}, \pi_{-s}; \xi_s) = 0$

A2. Fulfilled demand is bounded from below by a positive constant; that is

$$D_s \geq D_{\min} \in R_{++}^1$$

for all $s \in S$.

A3. Every observed demand function $h(\pi; \xi_s)$ is concave in $\pi$ for all feasible prices and all uncertainty vectors $\xi_s$.

A4. For any fixed $\bar{\pi}_{-s}$ and uncertainty factors $\xi_s$, every demand function is strictly monotonically decreasing in own price; that is

$$\int_{t_0}^{t_f} \exp(-t)[h_s(\pi_s^1, \bar{\pi}_{-s}; \xi_s) - h_s(\pi_s^2, \bar{\pi}_{-s}; \xi_s)](\pi_s^1 - \pi_s^1) dt < 0$$

for all $s \in S$ and all feasible prices $\pi_s^1 \neq \pi_s^2$.

A5. Every demand function is linear in uncertainty factor $\xi_s$, that is

$$\frac{\partial^2 h_s}{\partial \xi_s^2} = 0$$

for all $s \in S$.

A6. Each norm $\|\nabla_{\xi_s} h_s(\pi_s, \pi_{-s}; \xi_s)\|$ is monotone increasing in own price $\pi_s$ for all $s \in S$, all feasible prices $\pi$ and all uncertainty factors $\xi_s$.

Assumptions A1, A2, A3 and A4 were made by Perakis and Sood (2006). Assumptions A5 and A6 are unique to this exposition and may be motivated by considering the following commonly illustrative linear demand function with multiplicative uncertainty factor:

$$h_s(\pi_s, \pi_{-s}; \xi_s) = (\alpha_s - \beta_s \pi_s + \sum_{r \neq s} \gamma_{sr} \pi_r) \xi_s \quad (3)$$

where $\alpha_s, \beta_s, \gamma_{sr} \in R_{++}^1$ for all $s, r \in S$. By inspection (3) is a linear function of $\pi$ and a decreasing function of own price, so assumptions A3 and A4 are satisfied. It is also clear that (3) is a linear function of $\xi_s$, is a scalar in this specific case, and assumption A5 is satisfied. Moreover, for (3) we have

$$\nabla_{\xi_s} h_s(\pi_s, \pi_{-s}; \xi_s) = \alpha_s - \beta_s \pi_s + \sum_{r \neq s} \gamma_{sr} \pi_r$$

which is a scalar and leads to

$$\|\nabla_{\xi_s} h_s(\pi_s, \pi_{-s}; \xi_s)\| = \left| \alpha_s - \beta_s \pi_s + \sum_{r \neq s} \gamma_{sr} \pi_r \right| \quad (4)$$

By inspection (4) is an increasing function of own price and Assumption A6 is satisfied.

## 3 Best Response Problem

For any given but arbitrary $\pi_{-s}$, seller $s \in S$ seeks to solve the following infinite dimensional mathematical program

$$\max J_s(\pi_s, D_s) = \max \int_{t_0}^{t_f} \exp(-\rho t) \pi_s(t) D_s(t) dt \quad (5)$$

subject to

$$D_s \leq h_s(\pi_s, \pi_{-s}; \xi_s) \quad \forall \xi_s \in U_s \quad (6)$$

$$K_s = \int_{t_0}^{t_f} D_s(t)dt \quad (7)$$

$$\pi_s \geq \pi_{\min} \in R_+^1 \quad (8)$$

$$\pi_s \leq \pi_{\max} \in R_{++}^1 \quad (9)$$

$$D_s \geq D_{\min} \in R_{++}^1 \quad (10)$$

where $K_s$ is the initial endowment of inventory possessed by each seller $s \in S$. Note that

It is helpful to re-state problem (5) through (10) as

$$\max J_s(\pi_s, D_s) = \max \int_{t_0}^{t_f} \exp(-\rho t)\pi_s(t) D_s(t) dt \quad (11)$$
$$s.t. \quad (\pi_s, D_s) \in \Lambda_s(\pi_{-s}, \xi_s)$$

where

$$\Lambda_s(\pi_{-s}, \xi_s) \equiv \left\{ (\pi_s, D_s) : \pi_{\min} - \pi_s \leq 0, \pi_s - \pi_{\max} \leq 0, \int_{t_0}^{t_f} D_s(t)dt - K_s = 0, \right.$$
$$\left. D_{\min} - D_s \leq 0, D_s - h_s(\pi_s, \pi_{-s}; \xi_s) \leq 0 \right\}$$

is the strategy space for seller $s \in S$. We also define

$$\Lambda(\pi) = \{(\pi, D) : (\pi_s, D_s) \in \Lambda_s(\pi_{-s}, \xi_s) \quad \forall s \in S\} \quad (12)$$

The strategy space $\Lambda_s(\pi_{-s}, \xi_s)$ is convex for all $s \in S$, a result formalized in the following lemma:

**Lemma 1** *Seller's convex strategy space*. For each seller $s \in S$, take $h_s(\pi_s, \pi_{-s})$ to be concave in $\pi_s$ for all $\pi_{-s}$ (Assumption A3). Then the robust strategy space of each seller, $\Lambda_s(\pi_{-s}, \xi_s)$, is convex in $(\pi_s, D_s)$ for all $s \in S$. Furthermore, $\Lambda(\pi)$ is convex in $(\pi, D)$.

**Proof:** Note that all the constraint functions for seller $s \in S$ are convex functions. In particular, let us consider the constraint

$$g_s(D_s) = \int_{t_0}^{t_f} D_s(t)dt - K_s = 0$$

For arbitrary points $D_s^1$ and $D_s^2$ with $\lambda \in [0,1] \subset R_+^1$, we have

$$g_s[\lambda D_s^1(t) + (1-\lambda)D_s^2(t)] = \int_{t_0}^{t_f} [\lambda D_s^1(t) + (1-\lambda)D_s^2(t)]dt - K_s$$

$$= \lambda \int_{t_0}^{t_f} D_s^1(t)dt - \lambda K_s + (1-\lambda)\int_{t_0}^{t_f} D_s^2(t)dt - (1-\lambda)K_s$$

$$= \lambda \left\{ \int_{t_0}^{t_f} D_s^1(t)dt - K_s \right\} + (1-\lambda)\left\{ \int_{t_0}^{t_f} D_s^2(t)dt - K_s \right\}$$

$$= \lambda g_s(D_s^1) + (1-\lambda)g_s(D_s^2)$$

Thus, $g_s(D_s)$ is linear and hence convex. All constraint functions are convex for seller $s \in S$, therefore, $\Lambda_s(\pi_{-s}, \xi_s)$ is convex as a set. Moreover, the convexity of the $\Lambda_s(\pi_{-s}, \xi_s)$ for all $s \in S$ assures the convexity of $\Lambda(\pi)$. ∎

Furthermore, we know that, for the function spaces stipulated, the G-derivative of the criterion functional is

$$\begin{aligned}
&\delta J_s(\pi_s, D_s; \phi_\pi, \phi_D) \\
&= \lim_{\theta \to 0} \int_{t_0}^{t_f} \exp(-\rho t) \frac{(\pi_s + \theta \phi_\pi)(D_s + \theta \phi_D) - \pi_s D_s}{\theta} dt \\
&= \lim_{\theta \to 0} \int_{t_0}^{t_f} \exp(-\rho t) \frac{(\pi_s D_s + \theta^2 \phi_\pi \phi_D + \theta D_s \phi_\pi + \theta \pi_s \phi_D) - \pi_s D_s}{\theta} dt \\
&= \lim_{\theta \to 0} \int_{t_0}^{t_f} \exp(-\rho t) \frac{\theta \phi_\pi \phi_D + D_s \phi_\pi + \pi_s \phi_D}{1} dt \\
&= \int_{t_0}^{t_f} \exp(-\rho t)(D_s \phi_\pi + \pi_s \phi_D)dt \\
&= \int_{t_0}^{t_f} \left(\exp(-\rho t)D_s, \exp(-\rho t)\pi_s\right) \begin{pmatrix} \phi_\pi \\ \phi_D \end{pmatrix} dt
\end{aligned} \quad (13)$$

Of course

$$\delta J_s(\pi_s, D_s; \phi_\pi, \phi_D) = \int_{t_0}^{t_f} [\nabla_s J_s(\pi_s, D_s)\phi]dt \quad (14)$$

where

$$\phi = \begin{pmatrix} \phi_\pi \\ \phi_D \end{pmatrix}$$

Upon comparing (13) and (14) we see that

$$\nabla_s J_s(\pi_s, D_s) = \begin{pmatrix} \partial J_s/\partial \pi_s \\ \partial J_s/\partial D_s \end{pmatrix} = \begin{pmatrix} \exp(-\rho t)D_s \\ \exp(-\rho t)\pi_s \end{pmatrix} \quad (15)$$

Therefore, the first order condition, when $(\pi_s^*, D_s^*) \in \Lambda_s(\pi_{-s}, \xi_s)$ is a solution, takes the form

$$\delta J_s(\pi_s^*, D_s^*; \phi_\pi, \phi_D) \leq 0$$

for all feasible directions

$$\phi_\pi = \pi_s - \pi_s^*, \quad \phi_D = D_s - D_s^*$$

The above is easily re-stated in the more familiar form

$$\left\langle \nabla_s J_s(\pi_s^*, D_s^*), \begin{pmatrix} \pi_s - \pi_s^* \\ D_s - D_s^* \end{pmatrix} \right\rangle \leq 0$$

or equivalently as

$$\int_{t_0}^{t_f} \left[ \frac{\partial J_s(\pi_s^*, D_s^*)}{\partial \pi_s}(\pi_s - \pi_s^*) + \frac{\partial J_s(\pi_s^*, D_s^*)}{\partial D_s}(D_s - D_s^*) \right] dt \leq 0$$

$$\forall (\pi_s, D_s) \in \Lambda_s(\pi_{-s}, \xi_s)$$

In light of our knowledge of the gradient of the criterion $J_s$, this last quasi-variational inequality may be stated as

$$\int_{t_0}^{t_f} \exp(-\rho t)[D_s^* \cdot (\pi_s - \pi_s^*) + \pi_s^* \cdot (D_s - D_s^*)] dt \leq 0 \quad (16)$$

$$\forall (\pi_s, D_s) \in \Lambda_s(\pi_{-s}, \xi_s)$$

Statement (16) will not only be a necessary condition but also a sufficient condition if $J_s(\pi_s, D_s)$ is pseudo-concave on $\Lambda_s(\pi_{-s}, \xi_s)$. The pseudo-concavity of each seller's criterion will allow us to establish an equivalent variational inequality formulation of the generalized Nash equilibrium among sellers described by (11). To that end we state and prove the following result:

**Lemma 2** *Seller's criterion is pseudo-concave.* For each $s \in S$, the criterion $J_s(\pi_s, D_s)$ is pseudo-concave on $\Lambda_s(\pi_{-s}, \xi_s)$.

**Proof:** For the criterion $J_s(\pi_s, D_s)$ to be pseudo-concave on $\Lambda_s(\pi_{-s}, \xi_s)$, we must show that

$$\left\langle \nabla_s J_s(\pi_s^2, D_s^2), \begin{pmatrix} \pi_s^1 - \pi_s^2 \\ D_s^1 - D_s^2 \end{pmatrix} \right\rangle \leq 0 \Rightarrow J_s(\pi_s^2, D_s^2) \geq J_s(\pi_s^1, D_s^1) \quad (17)$$

for $(\pi_s^1, D_s^1), (\pi_s^2, D_s^2) \in \Lambda_s(\pi_{-s}, \xi_s)$. Property (17) is assured if $R_s(\pi_s, D_s) = \pi_s D_s$ is pseudo-concave at each instant of time $t \in [t_0, t_f]$. By Theorem 9 of Ferland (1972), we know that $Z_s = (-1) \cdot R_s(\pi_s, D_s)$ is pseudo-convex (and hence $R_s(\pi_s, D_s)$ is pseudo-concave), if the following matrix

$$M_s = \begin{pmatrix} 0 & \dfrac{\partial Z_s}{\partial \pi_s} & \dfrac{\partial Z_s}{\partial D_s} \\ \dfrac{\partial Z_s}{\partial \pi_s} & \dfrac{\partial^2 Z_s}{\partial \pi_s^2} & \dfrac{\partial^2 Z_s}{\partial \pi_s \partial D_s} \\ \dfrac{\partial Z_s}{\partial D_s} & \dfrac{\partial^2 Z_s}{\partial D_s \partial \pi_s} & \dfrac{\partial^2 Z_s}{\partial D_s^2} \end{pmatrix} = \begin{pmatrix} 0 & -D_s & -\pi_s \\ -D_s & 0 & -1 \\ -\pi_s & -1 & 0 \end{pmatrix}$$

has a determinant that is strictly negative. We note that

$$\det M_s = -2\pi_s D_s < 0$$

since each seller $s \in S$ is constrained to employ a strictly positive solution bounded away from the origin. ∎

We are now prepared to present and prove the following result:

**Theorem 1** *Quasi-variational inequality equivalent to best response.* The quasi-variational inequality (16) is equivalent to the infinite dimensional mathematical program (5), (6), (7), (8), (9) and (10).

**Proof:** We have already demonstrated that a solution of the best response problem will obey (16) in the discussion preceding Lemma 2. It remains to show that any solution of (16) is a solution of (5), (6), (7), (8), (9) and (10). Note that because of the pseudo-concavity of the best response criterion, as expressed by (17), the variational inequality (16) implies

$$\int_{t_0}^{t_f} \exp(-\rho t)\pi_s^*(t)D_s^*(t)dt = J_s(\pi_s^*, D_s^*) \geq J_s(\pi_s, D_s) \\ = \int_{t_0}^{t_f} \exp(-\rho t)\pi_s(t)D_s(t)dt \quad (18)$$

for all $(\pi_s, D_s) \in \Lambda_s(\pi_{-s}, \xi_s)$. Thus, a solution $(\pi_s^*, D_s^*) \in \Lambda_s(\pi_{-s}, \xi_s)$ of the variational inequality associated with seller $s \in S$ is a solution of the infinite dimensional mathematical program describing that seller's best response. ∎

Existence of a solution of the best response problem can be assured by the following result:

**Theorem 2** *Existence of a best response.* The best response problem (11) faced by seller $s \in S$ with the assumptions introduced in Section 2 has a solution.

**Proof:** We note that (16) may be stated as

$$\langle F_s(v_s^*), v_s - v_v^* \rangle \leq 0 \quad v_s, v_v^* \in \Lambda_s(v_{-s}, \xi_s) \subset V \equiv L^2[t_0, t_f] \times L^2[t_0, t_f] \quad (19)$$

when

$$F = \begin{pmatrix} \exp(-\rho t)D_s \\ \exp(-\rho t)\pi_s \end{pmatrix} \qquad v_s = \begin{pmatrix} \pi_s \\ D_s \end{pmatrix} \qquad v_{-s} = \pi_{-s}$$

By Theorem 2 of Browder (1968), we know that (19) has a solution because $\Lambda_s(v_{-s}, \xi_s)$ is a convex and compact subset of $V$ while $F$ is a continuous mapping of $\Lambda_s(v_{-s}, \xi_s)$ into $V^* = V$. ∎

We next make note of the following:

**Lemma 3** *Demand constraint (1) binds for some* $\xi_s \in U_s$. For every competitor $s \in S$, given her competitors' strategies $\pi_{-s}$, a solution $(\pi_s^*, D_s^*)$ of a best response problem having constraints expressed as $(\pi_s, D_s) \in \Lambda_r(\pi_{-r}, \xi_r)$ and whose criterion is monotonically increasing in own price satisfies $D_s^* = h_s(\pi_s^*, \pi_{-s}; \xi_s^*)$ for some $\xi_s^* \in U_s$.

**Proof:** This result is proven in Perakis and Sood (2006) for a discrete time formulation of the same problem; its adaptation to the continuous-time models considered herein is trivial and hence omitted for the sake of brevity. ∎

We are now in a position to observe that the seller's best response policy is a unique strategy:

**Lemma 4** *Seller's best policy is unique.* For each $s \in S$ and all $r \in S$ such that $r \neq s$ and $(\pi_r, D_r) \in \Lambda_r(\pi_{-r}, \xi_r)$, any solution of variational inequality (16) is a unique solution of the best response problem.

**Proof:** We recall the assumed monotonicity of demand for own arguments: for any $\xi_s \in U_s$

$$\int_{t_0}^{t_f} \exp(-\rho t)[h_s(\pi_s^1, \pi_{-s}^*; \xi_s) - h_s(\pi_s^2, \pi_{-s}^*; \xi_s)](\pi_s^1 - \pi_s^2)dt < 0, \forall (\pi_s^1, \pi_s^2) \quad (20)$$

For a given $\pi_s$ and $\pi_s^*$, we chose $\xi_s$ and $\xi_s^*$ as follows:

$$\xi_s = \arg\min_{\xi_s \in U_s} h_s(\pi_s, \pi_{-s}^*; \xi_s)$$

$$\xi_s^* = \arg\min_{\xi_s \in U_s} h_s(\pi_s^*, \pi_{-s}^*; \xi_s)$$

Thus, we have for any $\xi_s$ and $\xi_s^* \in U_s$

$$\int_{t_0}^{t_f} \exp(-\rho t)[h_s(\pi_s^*, \pi_{-s}^*; \xi_s) - h_s(\pi_s, \pi_{-s}^*; \xi_s)](\pi_s^* - \pi_s)dt < 0$$

$$\int_{t_0}^{t_f} \exp(-\rho t)[h_s(\pi_s, \pi_{-s}^*; \xi_s^*) - h_s(\pi_s^*, \pi_{-s}^*; \xi_s^*)](\pi_s - \pi_s^*)dt < 0$$

There are two cases:

1. When $\pi_s^* \geq \pi_s$,

$$h_s(\pi_s^*, \pi_{-s}^*; \xi_s^*) \leq h_s(\pi_s^*, \pi_{-s}^*; \xi_s) \leq h_s(\pi_s, \pi_{-s}^*; \xi_s)$$

then we have

$$\int_{t_0}^{t_f} \exp(-\rho t)[h_s(\pi_s^*, \pi_{-s}^*; \xi_s^*) - h_s(\pi_s, \pi_{-s}^*; \xi_s)](\pi_s^* - \pi_s)dt < 0$$

2. When $\pi_s \geq \pi_s^*$,

$$h_s(\pi_s, \pi_{-s}^*; \xi_s) \leq h_s(\pi_s, \pi_{-s}^*; \xi_s^*) \leq h_s(\pi_s^*, \pi_{-s}^*; \xi_s^*)$$

then we have

$$\int_{t_0}^{t_f} \exp(-\rho t)[h_s(\pi_s^*, \pi_{-s}^*; \xi_s^*) - h_s(\pi_s, \pi_{-s}^*; \xi_s)](\pi_s^* - \pi_s)dt < 0$$

In both cases,

$$\int_{t_0}^{t_f} \exp(-\rho t)[h_s(\pi_s^*, \pi_{-s}^*; \xi_s^*) - h_s(\pi_s, \pi_{-s}^*; \xi_s)](\pi_s^* - \pi_s)dt < 0 \quad (21)$$

From Lemma 3, we know $D_s^* = h_s(\pi_s^*, \pi_{-s}^*; \xi_s^*)$, thus (21) can be re-stated as

$$\int_{t_0}^{t_f} \exp(-\rho t)[D_s^* - h_s(\pi_s, \pi_{-s}^*; \xi_s)](\pi_s^* - \pi_s)dt < 0 \quad (22)$$

We next recall the variational inequality associated with the best response problem:

$$\int_{t_0}^{t_f} \exp(-\rho t)[D_s^* \cdot (\pi_s - \pi_s^*) + \pi_s^* \cdot (D_s - D_s^*)]dt < 0, \forall (\pi_s, D_s) \in \Lambda_s(\pi_{-s}, \xi_s) \quad (23)$$

Adding (22) and (23) yields

$$\int_{t_0}^{t_f} \exp(-\rho t)[h_s(\pi_s, \pi_{-s}^*; \xi_s)(\pi_s - \pi_s^*) + \pi_s^* \cdot (D_s - D_s^*)]dt < 0 \quad (24)$$
$$\forall (\pi_s, D_s) \in \Lambda_s(\pi_{-s}, \xi_s)$$

which may be re-stated as

$$\int_{t_0}^{t_f} \exp(-\rho t) Q_s(\pi_s, \pi_{-s}^*; \xi_s)dt < 0 \quad \forall (\pi_s, D_s) \in \Lambda_s(\pi_{-s}, \xi_s) \quad (25)$$

where

$$Q_s(\pi_s, \pi_{-s}^*; \xi_s) = h_s(\pi_s, \pi_{-s}^*; \xi_s)\pi_s + [D_s - h_s(\pi_s, \pi_{-s}^*; \xi_s)]\pi_s^* - \pi_s^* D_s^* \quad (26)$$

Define

$$\pi_s^0 = \sup_{\pi_s}\{\pi_s : D_s - h_s(\pi_s, \pi_{-s}^*; \xi_s) \leq 0\}$$

Upon taking $\pi_s = \pi_s^0$ and noting that Lemma 3 requires $D_s = h_s(\pi_s^0, \pi_{-s}^*; \xi_s)$, (26) becomes

$$Q_s(\pi_s^0, \pi_{-s}^*; \xi_s) = D_s \pi_s^0 - \pi_s^* D_s^* \geq D_s \pi_s - \pi_s^* D_s^* \quad (27)$$

In light of (27), inequality (25) yields

$$\int_{t_0}^{t_f} \exp(-\rho t) D_s \pi_s \, dt - \int_{t_0}^{t_f} \exp(-\rho t) D_s^* \pi_s^* \, dt \leq \int_{t_0}^{t_f} \exp(-\rho t) Q_s(\pi_s^0, \pi_{-s}^*; \xi_s) dt < 0$$

(28)

which establishes that a solution of variational inequality (16) is a unique global solution of the best response problem of competitor $s \in S$. ∎

## 4 Robust Best Response Problem

We next articulate a robust counterpart of the best response problem for seller $s \in S$. This is accomplished by adding to the demand constraints

$$D_s - h_s(\pi_s, \pi_{-s}; \xi_s) \leq 0 \quad \forall s \in S, \xi_s \in U_s$$

an additional term that represents a safety margin or so-called protection level. Based on the nominal values of uncertain parameters

$$\xi_s^0 \in (L^2[t_0, t_f])^m \quad \forall s \in S$$

we may introduce the uncertainty sets

$$U_s = \{\xi_s^0(t) + \tau\psi(t) : \|\psi(t)\| \leq 1 \quad \tau \in R_{++}^1\} \quad \forall s \in S \quad (29)$$

where

$$\|\psi(t)\|^2 = \int_{t_0}^{t_f} [\psi(t)]^T \psi(t) dt \quad (30)$$

and

$$\xi_s(t) \in (L^2[t_0, t_f])^m$$

$$\psi(t) \in (L^2[t_0, t_f])^m$$

The following result provides the robust constraint needed to articulate the robust best response problem for each seller $s \in S$:

**Lemma 5** *Robust counterpart of constraint (6).* For each seller $s \in S$, the set of constraints

$$\{D_s - h_s(\pi_s, \pi_{-s}; \xi_s) \leq 0 \quad \xi_s \in U_s\}$$

has the following robust counterpart:

$$D_s - h_s(\pi_s, \pi_{-s}; \xi_s^0) + \tau \left\| \nabla_{\xi_s} h_s(\pi_s, \pi_{-s}; \xi_s^0) \right\| \leq 0 \quad \forall \xi_s \in U_s$$

where the norm is taken in $(L^2[t_0, t_f])^m$.

**Proof:** In proving this lemma, we are influenced by the approach of Zhang (2007) for a similar but not identical result. Let us first introduce the following definition:

$$g_s(D_s, \pi_s, \pi_{-s}; \xi_s) \equiv D_s - h_s(\pi_s, \pi_{-s}; \xi_s) \leq 0 \quad \xi_s \in U_s\}$$

We seek

$$\max_{\xi_s \in U_s} g_s(D_s, \pi_s, \pi_{-s}; \xi_s) \leq 0$$

Using a first-order Taylor expansion around the nominal value $\xi_s^0$, we obtain

$$\max_{\xi_s \in U_s} g_s(D_s, \pi_s, \pi_{-s}; \xi_s) = g_s(D_s, \pi_s, \pi_{-s}; \xi_s^0)$$

$$+ \max_{\xi_s \in U_s} \int_{t_0}^{t_f} [\nabla_{\xi_s} g_s(D_s, \pi_s, \pi_{-s}; \xi_s^0)]^T (\xi_s - \xi_s^0) dt$$

$$= g_s(D_s, \pi_s, \pi_{-s}; \xi_s^0)$$

$$+ \tau \max_{\|\psi\| \leq 1} \int_{t_0}^{t_f} [\nabla_{\xi_s} g_s(D_s, \pi_s, \pi_{-s}; \xi_s^0)]^T \psi \, dt$$

(31)

where (31) is a maximization problem having $\psi$ as its vector of decision variables. Note that the above expansion is not a local approximation, because the function being expanded is linear in $\xi_s$ (Assumption A5). From the $L^2$-version of the well-known Schwarz's inequality[†], we have

$$\int_{t_0}^{t_f} \left| [\nabla_{\xi_s} g_s(D_s, \pi_s, \pi_{-s}; \xi_s^0)]^T \psi \right| dt$$

$$\leq \left( \int_{t_0}^{t_f} [\nabla_{\xi_s} g_s(D_s, \pi_s, \pi_{-s}; \xi_s^0)]^T [\nabla_{\xi_s} g_s(D_s, \pi_s, \pi_{-s}; \xi_s^0)] dt \right)^{\frac{1}{2}} \left( \int_{t_0}^{t_f} \psi^T \psi dt \right)^{\frac{1}{2}}$$

---

[†] See, for example, Rudin (1987).

Therefore,

$$\max_{\|\psi\| \leq 1} \int_{t_0}^{t_f} [\nabla_{\xi_s} g_s(D_s, \pi_s, \pi_{-s}; \xi_s^0)]^T \psi dt$$

$$\leq \max_{\|\psi\| \leq 1} \int_{t_0}^{t_f} \left| [\nabla_{\xi_s} g_s(D_s, \pi_s, \pi_{-s}; \xi_s^0)]^T \psi \right| dt$$

$$\leq \left( \int_{t_0}^{t_f} [\nabla_{\xi_s} g_s(D_s, \pi_s, \pi_{-s}; \xi_s^0)]^T [\nabla_{\xi_s} g_s(D_s, \pi_s, \pi_{-s}; \xi_s^0)] dt \right)^{\frac{1}{2}}$$

$$= \left\| \nabla_{\xi_s} g_s(D_s, \pi_s, \pi_{-s}; \xi_s^0) \right\| \qquad (32)$$

in light of (30). Consequently, we have

$$\max_{\xi_s \in U_s} g_s(D_s, \pi_s, \pi_{-s}; \xi_s) \leq g_s(D_s, \pi_s, \pi_{-s}; \xi_s^0) + \tau \left\| \nabla_{\xi_s} g_s(D_s, \pi_s, \pi_{-s}; \xi_s^0) \right\| \leq 0$$

Therefore, the robust counterpart of constraint (6) is

$$D_s - h_s(\pi_s, \pi_{-s}; \xi_s^0) + \tau \left\| \nabla_{\xi_s} h_s(\pi_s, \pi_{-s}; \xi_s^0) \right\| \leq 0$$

since

$$\left\| \nabla_{\xi_s} g_s(D_s, \pi_s, \pi_{-s}; \xi_s^0) \right\| = \left\| -\nabla_{\xi_s} h_s(\pi_s, \pi_{-s}; \xi_s^0) \right\| = \left\| \nabla_{\xi_s} h_s(\pi_s, \pi_{-s}; \xi_s^0) \right\|$$

∎

The robust best response problem may be stated as

$$\max J_s(\pi_s, D_s) = \max \int_{t_0}^{t_f} \exp(-\rho t) \pi_s(t) D_s(t) dt \qquad (33)$$

subject to

$$D_s - h_s(\pi_s, \pi_{-s}; \xi_s^0) + \tau \left\| \nabla_{\xi_s} h_s(\pi_s, \pi_{-s}; \xi_s^0) \right\| \leq 0 \qquad (34)$$

$$K_s = \int_{t_0}^{t_f} D_s(t) dt \qquad (35)$$

$$\pi_s \geq \pi_{\min} \in R_+^1 \qquad (36)$$

$$\pi_s \leq \pi_{\max} \in R_{++}^1 \qquad (37)$$

$$D_s \geq D_{\min} \in R_{++}^1 \qquad (38)$$

We would like the robust best response problem to be a generalized convex program, so that its necessary conditions are also sufficient and global solutions may be computed. In preparation for our analysis of the robust best response problem, we state and prove the following result:

Lemma 6 *Monotonicity of the robust constraint*. Let assumption A6 hold so that the norm

$$\left\|\nabla_{\xi_s} h_s(\pi_s, \pi_{-s}; \xi_s^0)\right\| \quad (39)$$

is monotone increasing in own price $\pi_s$. Then, for all $s \in S$, each robust constraint function

$$g_s^0(D_s, \pi_s, \pi_{-s}; \xi_s^0) \equiv D_s - h_s(\pi_s, \pi_{-s}; \xi_s^0) + \tau \left\|\nabla_{\xi_s} h_s(\pi_s, \pi_{-s}; \xi_s^0)\right\|$$

forming the lefthand-side of (34) is monotone increasing in $\pi_s$.

**Proof:** Clearly for arbitrary feasible $\pi_s^1$ and $\pi_s^2$, monotonicity of (39) requires

$$\left(\left\|\nabla_{\xi_s} h_s^1\right\| - \left\|\nabla_{\xi_s} h_s^2\right\|\right)(\pi_s^1 - \pi_s^2) \geq 0 \quad (40)$$

where we have used the abbreviations

$$\nabla_{\xi_s} h_s^k \equiv \nabla_{\xi_s} h_s(\pi_s^k, \pi_{-s}^k; \xi_s^0) \quad k = 1, 2$$

The monotone decreasing nature of $h_s(\pi_s, \pi_{-s}; \xi_s)$ with respect to $\pi_s$ is expressed as

$$\left(h_s^1 - h_s^2\right)(\pi_s^1 - \pi_s^2) \leq 0$$

or

$$\left(-h_s^1 + h_s^2\right)(\pi_s^1 - \pi_s^2) \geq 0 \quad (41)$$

Multiplying (40) by $\tau$ and adding the result to (41) gives

$$\left(-h_s^1 + h_s^2 + \tau\left\|\nabla_{\xi_s} h_s^1\right\| - \tau\left\|\nabla_{\xi_s} h_s^2\right\|\right)(\pi_s^1 - \pi_s^2) \geq 0$$

which may be arranged to reveal

$$\left[\left(D_s - h_s^1 + \tau\left\|\nabla_{\xi_s} h_s^1\right\|\right) - \left(D_s - h_s^2 + \tau\left\|\nabla_{\xi_s} h_s^2\right\|\right)\right](\pi_s^1 - \pi_s^2) \geq 0 \quad (42)$$

where we have added and subtracted $D_s$ from the leading term of (42). The desired result is proven. ∎

We next observe that

**Lemma 7** *Quasiconvexity of a monotone function of a single variable*. A monotone increasing function of a single variable is quasiconvex.

**Proof:** We consider $f : R^1 \to R^1$, $x^1, x^2 \in L^2[t_0, t_f]$

Assume $x^1 < x^2$, in light of the assumed monotonicity, we have

$$f(x^1) \leq f(x^2)$$

$$f(\lambda x^1 + (1-\lambda)x^2) \leq f(x^2)$$

when $\lambda \in [0,1]$. It is immediate that

$$f(\lambda x^1 + (1-\lambda)x^2) \leq \max[f(x^1), f(x^2)]$$

which is the defining property of a quasiconvex function. ∎

The results developed above lead to

**Lemma 8** *Quasiconvexity of robust constraints.* Under the assumptions of Lemma 6, each robust constraint function

$$g_s^0(D_s, \pi_s, \pi_{-s}; \xi_s^0) \equiv D_s - h_s(\pi_s, \pi_{-s}; \xi_s^0) + \tau \left\| \nabla_{\xi_s} h_s(\pi_s, \pi_{-s}; \xi_s^0) \right\|$$

forming the lefthand-side of (34) is quasiconvex in $\pi_s$ for all $s \in S$.

**Proof:** The result follows immediately from Lemma 6 and 7. ∎

Lemma 6 establishes that quasiconvexity of the robust demand constraints will occur so long as $D_{\min}$ is not too small relative to other parameters. However, as has already been remarked, the most important implication of Lemma 8 is that the robust best response problem given by (33) through (38) is a generalized convex program involving the maximization of a pseudo-concave objective subject to quasiconvex inequality constraints. In particular Lemma 8 allows us to establish the following result:

**Lemma 9** Feasible region of the robust best response problem is convex. For every seller $s \in S$, the set

$$\Lambda_s(\pi_{-s}, \xi_s^0) \equiv \left\{ (\pi_s, D_s) : \pi_{\min} - \pi_s \leq 0, \pi_s - \pi_{\max} \leq 0, \int_{t_0}^{t_f} D_s(t)dt - K_s = 0, \right.$$
$$\left. D_{\min} - D_s \leq 0, g_s^0(D_s, \pi_s, \pi_{-s}; \xi_s^0) \leq 0 \right\}$$

is convex.

**Proof:** It is well known (see Avriel, 1976, Theorem 6.1) that the level sets of a quasiconvex function are convex; hence

$$G_s^0(\pi_{-s};\xi_s^0) = \Big\{(\pi_s, D_s) : g_s^0(D_s, \pi_s, \pi_{-s};\xi_s^0) \equiv$$
$$D_s - h_s(\pi_s, \pi_{-s};\xi_s^0) + \tau \left\| \nabla_{\xi_s} h_s(\pi_s, \pi_{-s};\xi_s^0) \right\| \leq 0 \Big\}$$

is a convex set for every $s \in S$ by virtue of Lemma 8. Furthermore

$$E_s \equiv \Big\{(\pi_s, D_s) : \pi_{\min} - \pi_s \leq 0, \pi_s - \pi_{\max} \leq 0, \int_{t_0}^{t_f} D_s(t)dt - K_s = 0,$$
$$D_{\min} - D_s \leq 0 \Big\}$$

is convex by virtue of Lemma 8. Since

$$\Lambda_s(\pi_{-s}, \xi_s^0) = G_s^0(\pi_{-s}, \xi_s^0) \cap E_s$$

for $s \in S$, the desired result follows immediately. ∎

We are led by virtue of the above development to the following result:

**Theorem 3** *Quasi-variational inequality equivalent to robust best response.* The following infinite dimensional quasi-variational inequality is equivalent to the infinite dimensional mathematical program (33), (34), (35), (36), (37) and (38) that describes the robust best response of each seller $s \in S$:

$$\text{find} \ \ (\pi_s^*, D_s^*) \in \Lambda_s(\pi_{-s}, \xi_s^0) \ \text{such that}$$

$$\int_{t_0}^{t_f} \exp(-\rho t)[D_s^* \cdot (\pi_s - \pi_s^*) + \pi_s^* \cdot (D_s - D_s^*)]dt \leq 0 \quad \forall (\pi_s, D_s) \in \Lambda_s(\pi_{-s}, \xi_s^0)$$

(43)

where

$$\Lambda_s(\pi_{-s}, \xi_s^0) \equiv \Big\{(\pi_s, D_s) : \pi_{\min} - \pi_s \leq 0, \pi_s - \pi_{\max} \leq 0, \int_{t_0}^{t_f} D_s(t)dt - K_s = 0,$$
$$D_{\min} - D_s \leq 0, g_s^0(D_s, \pi_s, \pi_{-s};\xi_s^0) \leq 0 \Big\}$$

**Proof:** In light of Lemma 9, the arguments of Lemma 2 apply. ∎

Existence of a solution of the robust best response problem can be assured by the following result:

**Theorem 4** *Existence of a robust best response.* The robust best response problem (33), (34), (35), (36), (37) and (38) faced by seller $s \in S$ with the assumptions introduced in Section 2 has a solution.

**Proof:** The proof employs arguments identical to those used in the proof of Theorem 2.

# 5 Quasi-Variational Inequality Formulation of the Generalized Robust Nash Equilibrium

We are now interested in formulating the non-cooperative game among sellers described by the robust best response problem (33) through (38). We employ, as a solution concept for that game, the notion of a generalized robust Nash equilibrium. The relevant result is the following:

**Theorem 5** *Generalized Robust Nash equilibrium among sellers expressible as a quasi-variational inequality.* The generalized robust Nash equilibrium among sellers $s \in S$ that is described by best response problem (33), (34), (35), (36), (37) and (38), is equivalent to the following market quasi-variational inequality:

$$\text{find } (\pi^*, D^*) \in \Lambda(\pi^*, \xi^0) \text{ such that}$$

$$\sum_{s \in S} \int_{t_0}^{t_f} \exp(-\rho t)[D_s^* \cdot (\pi_s - \pi_s^*) + \pi_s^* \cdot (D_s - D_s^*)]dt \leq 0 \quad \forall (\pi, D) \in \Lambda(\pi^*, \xi^0)$$

(44)

where

$$\Lambda(\pi^*, \xi^0) = \prod_{s \in S} \Lambda_s(\pi_{-s}^*, \xi_s^0)$$

when the regularity conditions known as A1, A2, A3, A4, A5 and A6 hold.

**Proof:** First we note that summation of (43) over $s \in S$ yields (44); thus, the following concatenation of solutions of robust best response problems

$$(\pi^*, D^*) = (\pi_1^*, \cdots \pi_s^*, \cdots \pi_{|S|}^*; D_1^*, \cdots D_s^*, \cdots D_{|S|}^*)$$

is a solution of (44). It remains to show that any solution of (44) also solves the individual robust best response problems for all $s \in S$. To that end, we note that (44) may be re-stated as

$$\max \sum_{s \in S} \int_{t_0}^{t_f} \exp(-\rho t)(D_s^* \cdot \pi_s + \pi_s^* \cdot D_s)dt \leq 0 \quad (45)$$

subject to

$$(\pi_s, D_s) \in \Lambda_s(\pi_{-s}^*, \xi_s^0) \quad \forall s \in S \quad (46)$$

The mathematical program formed by (45) and (46) may be decomposed into seller-specific programs of the form

$$\max \int_{t_0}^{t_f} \exp(-\rho t)(D_s^* \cdot \pi_s + \pi_s^* \cdot D_s)dt \leq 0 \quad s.t. \quad (\pi_s, D_s) \in \Lambda_s(\pi_{-s}^*, \xi_s^0)$$

(47)

for each $s \in S$. By inspection of (47), we have that the 2-tuple $(\pi_s^*, D_s^*)$ must obey

$$\int_{t_0}^{t_f} \exp(-\rho t)(D_s^* \cdot \pi_s + \pi_s^* \cdot D_s) dt \le \int_{t_0}^{t_f} \exp(-\rho t)(D_s^* \cdot \pi_s^* + \pi_s^* \cdot D_s^*) dt$$

$$\forall (\pi_s, D_s) \in \Lambda_s(\pi_{-s}^*, \xi_s^0) \quad (48)$$

for each $s \in S$. From (48) we have for each $s \in S$ that

$$\int_{t_0}^{t_f} \exp(-\rho t)[D_s^* \cdot (\pi_s - \pi_s^*) + \pi_s^* \cdot (D_s - D_s^*)] dt \le 0 \quad \forall (\pi_s, D_s) \in \Lambda_s(\pi_{-s}^*, \xi_s^0)$$

(49)

Expression (49) is recognized as the variational inequality that is equivalent to the robust best response problem of seller $s \in S$. Thus, any solution of the market robust quasi-variational inequality provides solutions for the robust best response problems of all sellers. The proof is complete. ∎

We also want to establish that the market equilibrium described by (44) exists. That result is provided by the following theorem:

**Theorem 6** *Existence of a generalized robust Nash equilibrium.* The market quasi-variational inequality (44), for the assumptions introduced in Section 2 has a solution. As a result, a generalized robust Nash equilibrium exists.

**Proof:** We note that (44) may be stated as

$$\langle F(v^*), v - v^* \rangle \le 0 \quad v, v^* \in \Lambda(\pi^*, \xi^0) \subset V \equiv (L^2[t_0, t_f])^{2|S|} \quad (50)$$

when

$$F = \begin{pmatrix} \exp(-\rho t) D \\ \exp(-\rho t) \pi \end{pmatrix} \quad v = \begin{pmatrix} \pi \\ D \end{pmatrix}$$

By Theorem 2 of Browder (1968), we know that (50) has a solution because $\Lambda(\pi^*, \xi^0)$ is a convex and compact subset of $V$ while $F$ is a continuous mapping of $\Lambda(\pi^*, \xi^0)$ into $V^* = V$. In light of Theorem 5, a generalized robust Nash equilibrium exists. ∎

## 6 Solving the Market Quasi-Variational Inequality

In this section, we wish to devise an algorithm for solving the generalized robust Nash equilibrium for competitive dynamic pricing with fixed inventories. Since the quasi-variational inequality (44) is our representation of the generalized robust Nash equilibrium, we need a numerical scheme for solving a continuous-time quasi-variational inequality. Efficient convergent methods for solving general quasi-variational inequalities are not available. Therefore, we restrict ourselves, for reasons that will be made clear shortly, to the following class of demand constraint functions for every $s \in S$:

$$g_s^0(D_s, \pi_s, \pi_{-s}; \xi_s) \equiv D_s - g_s^1(\pi_s; \xi_s^0) - g_s^2(\pi_{-s}; \xi_s^0) \qquad (51)$$

where

$g_s^1(\pi_s; \xi_s^0)$ =an instantaneously monotone decreasing function of own price $\pi_s$

$g_s^2(\pi_{-s}; \xi_s^0)$ =an instantaneously monotone increasing function of non-own prices $\pi_{-s}$

In that $g_s^1(\pi_s; \xi_s^0)$ is related to own-demand and $g_s^2(\pi_{-s}; \xi_s^0)$ to non-own demand, the monotone properties given to them above are behaviorally realistic.

## 6.1 Variational Inequality Solutions

We will make use a trivial extension of a result first given by Harker (1991). In abstract form, the result pertains to the following quasi-variational inequality:

$$\text{find } (x^*, y^*) \in \Lambda(x^*) \text{ such that}$$

$$\langle F_1(x^*, y^*), x - x^* \rangle + \langle F_2(x^*, y^*), y - y^* \rangle \geq 0 \quad \forall (x, y) \in \Lambda(x^*) \qquad (52)$$

where

$$F_1 : (L^2[t_0, t_f])^{n_1} \times (L^2[t_0, t_f])^{n_1} \to R^{n_1}$$

$$F_2 : (L^2[t_0, t_f])^{n_1} \times (L^2[t_0, t_f])^{n_2} \to R^{n_2}$$

$$x, x^* \in (L^2[t_0, t_f])^{n_1}$$

$$y, y^* \in (L^2[t_0, t_f])^{n_2}$$

The result of interest is the following:

**Theorem 7** *Variational inequality solutions of quasi-variational inequalities.* If there is a set $\Gamma$ such that

$$(i) \Lambda(x) \subseteq \Gamma \quad \forall x \in \Gamma$$

$$(ii) \ x \in \Lambda(x) \quad \forall x \in \Gamma$$

then any solution of the variational inequality: find $(x^*, y^*) \in \Gamma$ such that

$$\langle F_1(x^*, y^*), x - x^* \rangle + \langle F_2(x^*, y^*), y - y^* \rangle \geq 0 \quad \forall (x, y) \in \Gamma$$

is a solution of the quasi-variational inequality (52).

**Proof:** This is a trivial extension of Theorem 3 of Harker (1991), and the proof given there is easily adapted to prove the above result. ∎

Now, following Harker (1991), we let $\Gamma$ have the form

$$\Gamma = \{(x, y): y \in \Psi, g(x, y) \leq 0\} \quad (53)$$

with

$$g(x, y) = g^1(x, y) + g^2(x, y) \quad (54)$$

If we define

$$\Lambda(x^*) = \{(x, y): y \in \Psi, g^1(x, y) + g^2(x^*, y) \leq 0\} \quad (55)$$

then conditions (*i*) and (*ii*) of Theorem 7 are satisfied. Furthermore the abstract sets $\Gamma$ and $\Lambda(x)$ have a structure that mimics the constraints of the fixed inventory dynamic pricing problem of interest here. Thus, Theorem 7 assures that any solution of the following variational inequality will be a solution of the quasi-variational inequality (44):

$$\text{find } (\pi^*, D^*) \in \Gamma \text{ such that}$$

$$\sum_{s \in S} \int_{t_0}^{t_f} \exp(-\rho t)[D_s^* \cdot (\pi_s - \pi_s^*) + \pi_s^* \cdot (D_s - D_s^*)]dt \leq 0 \quad \forall (\pi, D) \in \Gamma \quad (56)$$

where

$$\Gamma = \Big\{(\pi, D): \pi_{\min} - \pi_s \leq 0, \pi_s - \pi_{\max} \leq 0, \int_{t_0}^{t_f} D_s(t)dt - K_s = 0, D_{\min} - D_s \leq 0,$$
$$D_s - g_s^1(\pi_s; \xi_s^0) - g_s^2(\pi_{-s}; \xi_s^0) \leq 0 \quad \forall s \in S\Big\}$$

and

$$\Lambda(\pi^*) = \Big\{(\pi, D): \pi_{\min} - \pi_s \leq 0, \pi_s - \pi_{\max} \leq 0, \int_{t_0}^{t_f} D_s(t)dt - K_s = 0, D_{\min} - D_s \leq 0,$$
$$D_s - g_s^1(\pi_s; \xi_s^0) - g_s^2(\pi_{-s}^*; \xi_s^0) \leq 0 \quad \forall s \in S\Big\}$$

The preceding discussion constitutes a constructive result of the following:

**Theorem 8** *VI that solves the QVI.* Any solution of the variational inequality (56) will solve the quasi-variational inequality (44).

## 6.2 Equivalent Fixed Point Problem

It is possible to solve (56) using a simple but effective fixed point algorithm. To understand the fixed point algorithm is helpful to consider the following abstract variational inequality: find $u^* \in U$ such that

$$\langle F(u^*,t), u - u^* \rangle \leq 0 \quad \forall u \in U \qquad (57)$$

where

$$F : (L^2[t_0, t_f])^m \times R^1_+ \to (H^1[t_0, t_f])^m$$

$$U \subseteq (L^2[t_0, t_f])^m$$

We have the following result:

**Theorem 9** *Fixed point formulation of infinite dimensional variational inequality.* If $F(u,t)$ is differentiable and convex with respect to $u$ while $U$ is convex, the variational inequality (56) is equivalent to the following fixed point problem:

$$u^* = P_U[u - \alpha F(u^*, t)] \qquad (58)$$

where $P_U[\cdot]$ is the minimum norm projection onto $U$ and $\alpha \in R^1_{++}$ is an arbitrary positive constant.

**Proof:** The fixed point problem under consideration requires that

$$u^* = \arg\min_u \left\{ \frac{1}{2} \| u^* - \alpha F(u^*, t) - u \|^2 : v \in U \right\} \qquad (59)$$

where $\alpha \in R^1_{++}$ is any positive real scalar. That is, to solve the fixed point problem, we seek the solution of

$$\min_u Z = \int_{t_0}^{t_f} \frac{1}{2} [u^* - \alpha F(u^*, t) - u]^T [u^* - \alpha F(u^*, t) - u] dt \quad s.t. \quad u \in U$$

where $u^*$ is treated as fixed for the purpose of projection. A valid first order condition that is both necessary and sufficient due to convexity of $F(u,t)$ and $U$ is the following: find $u^0 \in U$ such that

$$\int_{t_0}^{t_f} \left\{ \nabla_v \frac{1}{2} [u^* - \alpha F(u^*, t) - u^0]^2 \right\} (u - u^0) dt \geq 0 \qquad \forall u \in U \qquad (60)$$

By (59) we know that $u^* = u^0$, so (60) reduces to

$$\int_{t_0}^{t_f} (-1)\alpha F(u^*,t)(u-u^*)dt \geq 0 \qquad \forall u \in U \qquad (61)$$

which is evidently equivalent to (59). This completes the proof. ∎

## 6.3 A Fixed Point Algorithm

Naturally there is a fixed point algorithm associated with the above fixed point formulation and expressed as the following iterative scheme:

$$u^{k+1} = P_U[u^k - \alpha F(u^k,t)] \qquad (62)$$

The positive scalar $\alpha \in R^1_{++}$ may be chosen empirically to assist convergence and may even be changed as the algorithm progresses. We give the following detailed statement of the fixed point algorithm:

*Fixed Point Algorithm*

**Step 0.** *Initialization.* Identify an initial feasible solution $u^0$ such that $u^0 \in U$ and set $k = 0$.

**Step 1.** *Solve optimal control problem.* Solve the following optimal control problem:

$$\min_v J^k(v) = \int_{t_0}^{t_f} \frac{1}{2}[u^k - \alpha F(u^k,t) - v]^2 dt \geq 0 \qquad s.t. \quad v \in U \qquad (63)$$

Call the solution $u^{k+1}$.

**Step 2**. *Stopping test.* If $\left\| u^{k+1} - u^k \right\| \leq \varepsilon_1 \|$ where $\varepsilon_1 \in R^1_{++}$ is a preset tolerance, stop and declare $u^* \approx u^{k+1}$. Otherwise set $k = k+1$ and go to Step 1.

The convergence of the above algorithm is guaranteed under certain conditions, which is shown in Friesz (2010).

## 6.4 Solving the Sub-Problems

It is important to realize that the fixed point algorithm (62) can be carried out in continuous time provided we employ a continuous-time representation of the solution of each subproblem (63). The subproblems are linearly constrained quadratic programs and may be solved in a variety of ways. For our purposes in this paper, we use a discrete time approximation of each subproblem that is solved using MINOS. A fifth order polynomial in time is then fit to the discrete time solution of each subproblem, and the next fixed point iteration is carried out in continuous time. Many other schemes may be invoked for solving the subproblem (63).

# 7 Pricing Rules and Dual Time Scales

Generally price changes take time to implement. Furthermore, in many circumstances, it is desirable or even mandatory to announce changes in pricing policy in advance. These considerations have the effect of introducing a second time scale, with the consequence that pricing involves explicit time shifts. An example is

$$\pi_s(t+\Delta) = \pi_s(t) + \sigma_s[h_s(\pi_s(t), \pi_{-s}(t); \xi_s) - D_s(t)] \quad \forall s \in S \quad \text{(abc)}$$

where $\Delta \in R^1_{++}$ is a positive, real constant representing the time needed to effect a price change. Furthermore, $\sigma_s \in R^1_{++}$ is a positive, real constant expressing the strength of the pricing response of seller $s \in S$ to changes in the excess of observed demand relative to realized demand.

Another pricing rule is

$$\pi_s(t+\Delta) - \pi_s(t) \geq 0 \quad \forall s \in S \quad \text{(???)}$$

which expresses, of course, monotonic pricing behavior, where again we assume that there is an intrinsic time needed to effect a price change. Still another pricing rule is

$$\pi_s(t+\Delta) \leq \frac{1}{t-t_0} \sum_{k \in S} \int_{t_0}^{t} \pi_k(\tau) d\tau \quad \forall s \in S \quad \text{(???)}$$

which stipulates that prices are set at or below the market's moving average price to reflect each seller's learning process relative to its competitors' prices, as well as the intrinsic time delay for price changes.

The above are but three instances of pricing rules, and many more such rules may be conjectured and imposed. It is, of course, important to know how the theory and computational approaches presented previously are impacted by the presence of such pricing rules. To that end, let us concentrate on constraints (abc), which are linear in $\pi_s(t)$ and $\pi_s(t+\Delta)$ for given $\pi_{-s}(t)$ when $h_s(\pi_s(t), \pi_{-s}(t); \xi_s)$ is linear. In that case, we may be certain that the strategy set of each seller is convex in both shifted and unshifted prices, as well as in realized demand. That is, the strategy set $\Lambda_s(\pi_{-s}, \xi_s)$, is convex in $(\pi_s, \pi_s^\Delta, D_s)$ for each $s \in S$, where we have used the convenient shorthand

$$\pi_s^\Delta(t) \equiv \pi_s(t+\Delta) \quad \forall s \in S \quad \text{(???)}$$

In the event that the observed demand function is linear in own price, constraint (abc) is a linear equality constraint for the best response problem of seller $s \in S$. As such it does not alter the linear-quadratic nature of subproblem (63). In particular, when computing, one may use an implicit fixed point perspective, which may be expressed as

$$\pi_s^{\Delta,k+1}(t+\Delta) \cong \pi_s^k(t+\Delta) \quad \forall s \in S \quad \text{(???)}$$

where this last expression is understood to mean that in iteration $k+1$ one uses the approximate price function found in iteration $k$.

# 8 Numerical Examples

In this section, we will present three numerical examples to show our results. The first example is used to show the general results of the equilibrium prices and profit. The second example will be used to compare the results of our robust pricing and non-robust pricing policy. The third example will illustrate the sensitivities analysis of demand parameters and robust parameters. We consider a linear demand function with multiplicative uncertainty

$$h_s(\pi_s, \pi_{-s}; \xi_s) = (\alpha_s - \beta_s \pi_s + \sum_{r \neq s} \gamma_{sr} \pi_r) \xi_s$$

where $\alpha_s, \beta_s, \gamma_{sr}$ are nonnegative constants and $\xi_s(t)$ is uncertain factor with nominal value $\xi_s^0(t)$ and magnitude $\tau$. The realization of $\xi_s(t)$ is randomly generated from the uncertainty set $U_s$, which is defined in (29).

The time horizon is $[t_0, t_f] = [1, 10]$. We consider two firms. The initial endowment of inventory for the two firms is $K_1 = 2500$ and $K_2 = 3000$, respectively. Assume the nominal value of the uncertainty factor $\xi_s^0(t) = 3 + 0.1t$, $s = \{1,2\}$, $t \in [t_0, t_f]$ and the magnitude $\tau = 0.8$.

**8.1 General Results of the Robust Pricing Policy**

*8.1.1 Two-firms example with identical demand function*

First we consider two competing sellers with the same linear demand functions. The parameters values are chosen as following:

$$\alpha_s(t) = 3000 \quad s = \{1,2\} \quad t \in [t_0, t_f]$$
$$\beta_s(t) = 180 - 4t \quad s = \{1,2\} \quad t \in [t_0, t_f]$$
$$\gamma_s(t) = 36 - 2t \quad s = \{1,2\} \quad t \in [t_0, t_f]$$

The results of the robust problem are shown in Figure 1, 2, 3 and Table 1. Figure 1 shows that the prices for both sellers are increasing with time. Seller 1's price is slightly higher than seller 2's price. Assume the uncertainty factor $\xi_s(t)$ follows uniform distribution in the uncertainty set, we run 10,000 simulations of the uncertainty factor, and then we can obtain the realized demand and compute the realized profit for each seller for each simulation. From Figure 2, we can see that the realized demand curve for each seller is quite stable. The value of realized demand for seller 1 and seller 2 at every instant time over the entire time horizon is around 250 and 300, respectively. Figure 3 shows the inventory change for each seller. The inventories of both sellers decrease to zero in the end of the time horizon. Before the end of the time horizon, seller 1' inventory is always less than seller 2' inventory. Combining Figure 1 and 3, we find that seller with less inventory (seller 1) sets higher price than seller with more inventory (seller 2), but seller 2 still obtains higher profit than seller 1, which is shown in Table 1.

*8.1.2 Two-firms example with unequal demand function*

We consider two firms with the following different parameters for the linear demand function:

$$\alpha_1(t) = 2500 \quad \alpha_2(t) = 3000 \quad t \in [t_0, t_f]$$
$$\beta_1(t) = 180 - 4t \quad \beta_2(t) = 170 - 4t \quad t \in [t_0, t_f]$$
$$\gamma_1(t) = 36 - 2t \quad \gamma_2(t) = 34 - 2t \quad t \in [t_0, t_f]$$

Figure 4 shows that the price for each seller is increasing with time. And seller 2's price is higher than seller 1. The realized demand and inventory change for each seller in this example are the same with those in Section 8.1.1. For brief, we omit the corresponding figures. The profits for both sellers are shown in Table 2.

## 8.2 Effect of Robustness on Price and Profit

In this section, we consider the effect of robustness on the equilibrium prices and profit. We consider two sellers with the following values for the linear demand function:

$$\alpha_1(t) = 2500 \quad \alpha_2(t) = 3500 \quad t \in [t_0, t_f]$$
$$\beta_1(t) = 175 - 4t \quad \beta_2(t) = 170 - 4t \quad t \in [t_0, t_f]$$
$$\gamma_1(t) = 35 - 2t \quad \gamma_2(t) = 34 - 2t \quad t \in [t_0, t_f]$$

In order to illustrate the effect of robustness, we design the following four scenarios, which are shown in Table 3, where `N' means the nominal pricing policy that the seller ignores the uncertainty in the demand function and obtains the optimal prices by naively assuming $\xi_s(t) = \xi_t^0(t) = 3 + 0.1t$. And `R' refers to our robust pricing policy.

Assume each seller believes that his competitor are using the same pricing policy (Bertsimas and Perakis, 2006), we can first obtain the equilibrium prices for each seller in Scenario 1 and 4, which are shown in Figure 5 and 6. Prices for each seller in Scenario 2 and 3 can be easily generated according to Figure 5 and 6. Seller 2' price is higher than seller 1 in both Scenario 1 and 4. We find that the prices for both sellers when both sellers are using the nominal pricing policy are higher than those when both sellers are using the robust pricing policy.

Here we assume $\xi_s(t)$ follows beta distribution and show the results of two different beta distributions, beta (1, 3) and beta (1, 1) (uniform distribution). We run 10,000 simulations for each distribution. Table 4 shows the range, the average value and the standard deviation of the realized profit for both sellers for the four scenarios with beta (1, 1) distribution. We can see that in Scenario 1 (Both sellers are using the nominal pricing policy), the range of profit is [42,271, 49,533] for seller 1 and [66,247, 79,003] for seller 2. The average profit for seller 1 is 46,730 with a standard deviation 1,145, and for seller 2 is 74,528 with a standard deviation 1,794. In Scenario 2 (Seller 2 changes to use the robust pricing policy), the profit for seller 1 is dramatically less than that in Scenario1, with average profit 22,510 and standard deviation 932. It implies that when his competitor adopts the robust pricing policy, a seller with nominal pricing policy performs much worse than that when he and his competitor both adopt the nominal pricing policy. In Scenario 2, although the average profit of seller 2 is lower than that in Scenario 1, seller 2 can always make a profit at 73, 885 no matter what the uncertainty is. Comparing the results of Scenario 3 to Scenario 1, we can draw the same conclusion. In Scenario 4, both sellers are using the robust pricing policy; seller 1 and 2 can obtain profit 45,928 and 73,885, respectively, both with standard deviation 0. Our robust pricing policy can guarantee sellers to achieve a profit with standard deviation 0, and improve much better the worst-case performance compared to the nominal pricing policy.

Table 5 shows the results with beta (1, 3) distribution. Comparing Table 5 and 4, we can obtain some new interesting conclusions: (1) The nominal pricing policy with beta (1, 3) performs worse than that with beta (1, 1), which implies that the performance of the nominal pricing policy depends on the distribution of the uncertainty factor; (2) No matter for beta (1, 1) or beta (1, 3), the robust pricing policy yields the same profit, which means the performance of the robust pricing policy does not depend on the distribution of the uncertainty factor. This makes our robust pricing policy more attractive since sellers typically do not know the exact distribution of the uncertainty; (3) The robust pricing policy achieves higher average profit than nominal pricing policy with beta (1, 3) distribution, no matter which pricing policies his competitor is using. For example, in Scenario 3, seller 2's average profit is 42,075, however, in Scenario 4, it increases to 73,855.

In order to show the results more detailed and accurate, we chose seller 1 arbitrarily and draw his profit histogram in Scenario 1, 2 and 3 with beta (1, 3) distribution, which is shown in Figure 7.

**8.3 Sensitivities Analysis**

*8.3.1 Sensitivities analysis of demand parameters*

In this section, we study the sensitivities analysis of demand parameters on price and profit. Assume the benchmark of the demand parameters take the following value:

$$\alpha_1(t) = 2500 \quad \alpha_2(t) = 3500 \quad t \in [t_0, t_f]$$
$$\beta_1(t) = 170 - 4t \quad \beta_2(t) = 170 - 4t \quad t \in [t_0, t_f]$$
$$\gamma_1(t) = 40 - 2t \quad \gamma_2(t) = 40 - 2t \quad t \in [t_0, t_f]$$

First, we vary $\alpha_1$ and see its effect on the price and profit of both sellers. Figure 8 and 9 show that both sellers' prices are increasing with $\alpha_1$, but seller 2's price increases less than seller 1. Figure 10 shows that the profits of both sellers are increasing with $\alpha_1$, but seller 2 increases less than seller 1, which means that seller 2 is much less sensitive to the changes of $\alpha_1$ than seller 1.

Second, we vary $\beta_1$ and see its effect on the price and profit of both sellers. Figure 11 and 12 shows that prices for seller 1 and 2 are decreasing with $\beta_1$, but seller 2 decreases less than seller 1. Figure 13 shows that the profits of seller 1 and 2 are also decreasing with $\beta_1$, but seller 2 decreases less than seller 1, which can be also explained that seller 2 is much less sensitive to the changes of $\beta_1$ than seller 1.

Finally, we vary $\gamma_2$ and see its effect on the price and profit of both sellers, which are shown in Figure 14-16. We can obtain the similar results as the analysis of $\alpha_1$.

*8.3.2 Sensitivities analysis of robust parameter*

In this section, we assume demand parameters for both sellers take the same value as in Section 8.2. We focus on the sensitivities analysis of robust parameter on profit. We will measure the profit when the actual uncertainty factor $\xi_s(t)$ can take any value from the uncertainty set while the seller employs a policy that assumes the uncertainty factor only takes value from a smaller set. We can do that by varying the magnitude $\bar{\tau}$ from 0 (nominal pricing policy) to 0.8 (very robust pricing policy), while the true magnitude $\tau = 0.8$. We obtain the corresponding prices according to $\bar{\tau}$, and calculate the realized profit. Assume $\xi_s(t)$ follows beta (1, 3) distribution. We consider the following three cases and run 10,000 simulations for each case:

I Seller 2 varies his magnitude while seller 1 adopts the nominal pricing policy

II Seller 2 varies his magnitude while seller 1 adopts the very robust pricing policy

III Seller 2 and seller 1 vary their magnitude simultaneously

Figure 17, 18 and 19 show the profit histogram for seller 2 in Case I, II and III, respectively. The three figures are all ordered by row-wise from left to right, corresponding to the magnitude $\bar{\tau}$ from 0 to 0.8. In Case I and III, we note that the standard deviation of the profit is reducing as seller 2 increases the robust parameter. The first six graphs of both Figure 17 and 19 show that the average profit is increasing with decreasing standard deviation, however, the last three graphs show that the average profit is decreasing with decreasing standard deviation, which means that upon a certain level of magnitude, if a seller wants to reduce the standard deviation of his profit, he needs to sacrifice a slight of his average profit. Thus, the seller has to tradeoff between average profit and standard deviation.

Figure 18 shows that the average profit of seller 2 first increases, and then decreases as seller 2 increases his robust parameter. However, before seller 2's magnitude $\bar{\tau}$ approaches to 0.8, the average profit of seller 2 in Case II is markedly lower than that in Case I and III, and the standard deviation in Case II reduces slower than that in Case I and III. This is due to that in Case II, seller 2's competitor adopts the very robust pricing policy, which makes seller 2 worse off than other cases until seller 2 adopts the very robust pricing policy as well.

# 8 Conclusions

We model a dynamic pricing problem in an oligopolistic market with demand uncertainty from a continuous-time perspective. Each seller with a fixed inventory competes with each other by setting prices. We use a multiplicative demand function to model the uncertain demand, and apply the ideas of robust optimization to deal with uncertainty. The uncertain demand constraint is replaced by a robust counterpart; hence the robust best response problem becomes a deterministic optimization problem. We have shown that the sellers' robust best response problem is a generalized Nash equilibrium problem (GNEP), and a quasi-variational inequality (QVI) is employed to formulate the GNEP. We prove that the QVI is equivalent to the GNEP, and a generalized robust Nash equilibrium exists among sellers. However, since there is no efficient convergent algorithm to solve QVI, we build a variational inequality (VI) and prove that any solution of the VI is a solution of the QVI, then use a fixed point algorithm to solve the VI. In the end, three numerical examples are presented to show the general results of our robust pricing policy, the effect of robustness on prices and profit and the sensitivities analysis of demand and robust parameters.

# References

[1] Aghassi, M., Bertsimas, D. (2006). Robust game theory. *Mathematical Programming, Ser. B*, 107, 231-273.

[2] Avriel, M. (1976). *Nonlinear Programming: Analysis and Methods*, Prentice Hall.

[3] Ben-Tal, A., Nemirovski, A. (1998). Robust convex optimization. *Mathematics of Operations Research*, 23, 769-805.

[4] Ben-Tal, A., Nemirovski, A. (1999). Robust solutions to uncertain linear programs. *Operations Research Letters*, 25, 1-13.

[5] Bertsimas, D., Perakis, G. (2006). Dynamic pricing: a learning approach. In S. Lawphongpanich, D. W. Hearn, & M. J. Smith (Eds.), *Mathematical and Computational Models for Congestion Charging* (pp. 45-80). New York: Springer.

[6] Bertsimas, D., Thiele, A. (2006). A robust optimization approach to inventory theory. Operations Research, 54(1), 150-168.

[7] Bitran, G., Caldentey, R. (2003). An overview of pricing models for revenue management. Manufacturing & Service Operations Management, 5, 203-229.

[8] Browder, F. E. (1968). The Fixed Point Theory of Multi-Valued Mappings in Topological Vector Spaces, Mathematische Annalen, 177, 283--301.

[9] Chen, X., Simchi-Levi, D. (2004). Coordinating inventory control and pricing strategies with random demand and fixed ordering cost: The finite horizon case. Operations Research, 52(6), 887-896.

[10] Elmaghraby, W., Keskinocak, P. (2003). Dynamic pricing in the presence of inventory considerations: Research overview, current practices, and future directions. Management Science, 49, 1287-1309.

[11] Feng, Y., Xiao, B., (2000). A continuous-time yield management model with multiple prices and reversible price changes. Management Science, 46 (5), 644-657.

[12] Ferland, J. A. (1972). Mathematical Programming with Quasi-Convex Objective Functions, Mathematical Programming, 3, 296-301.

[13] Friesz, T. L. (2010). *Dynamic optimization and differential games*. New York: Springer.

[14] Gallego, G., Hu, M. (2009). Dynamic pricing of perishable assets under competition. Columbia University, working paper.

[15] Gallego, G., van Ryzin, G. (1994). Optimal dynamic pricing of inventories with stochastic demand over finite horizons. Management Science, 40 (8), 999-1020.

[16] Gallego, G., Krishnamoorthy, S., Phillips, R. (2006). Dynamic revenue management games with forward and spot markets. Journal of Revenue and Pricing Management, 5(1), 10-31.

[17] Granot, D., Granot, F., Mantin, B. (2011). A dynamic pricing model under duopoly competition. University of British Columbia, Working paper.

[18] Guignard, M. (1969). Generalized Kuhn-Tucker Conditions for Mathematical Programming Problems in a Banach Space. SIAM Journal on Control, 7(2), 232-241.

[19] Harker, P. T. (1991). Generalized Nash Games and Quasi-Variational Inequalities. European Journal of Operational Research, 54, 81-94.

[20] Kachani, S., Perakis, G., Simon, C. (2004). A transient model for joint pricing and demand learning under competition. in: Presented in the Fourth Annual INFORMS Revenue Management and Pricing Section Conference, MIT, Cambridge.

[21] Kwon, C., Friesz, T.L., Mookherjee, R., Yao, T., Feng, B. (2009). Non-cooperative competition among revenue maximizing service providers with demand learning. European Journal of Operational Research, 197, 981-996.

[22] Leung, S. C. H., Tsang, S. O. S., Ng, W. L., Wu, Y. (2007). A robust optimization model for multi-site production planning problem in an uncertain environment. European Journal of Operational Research, 181, 224-238.

[23] Levin, Y., McGill, J., Nediak, M. (2009). Dynamic pricing in the presence of strategic consumers and oligopolistic competition. Management Science, 55 (1), 32-46.

[24] Levin, Y., McGill, J., Nediak, M. (2010). Optimal dynamic pricing of perishable items by a monopolist facing strategic consumers. Production and Operations Management, 19(1), 40-60.

[25] Lin, K.Y., Sibdari, S.Y. (2009). Dynamic price competition with discrete customer choices. European Journal of Operational Research, 197, 969-980.

[26] McGill, J. I., and van Ryzin, G. J. (1999). Revenue management: Research overview and prospects. Transportation Science, 33(2), 233-256.

[27] Pang, J. S., Stewart, D. (2008). Differential Variational Inequalities. Mathematical Programming, Ser. A 113, 345-424.

[28] Perakis, G., Sood, A. (2006). Competitive Multi-period Pricing for Perishable Products: A Robust Optimization Approach. Mathematical Programming, Ser. B 107, 295-335.

[29] Rudin, W. (1986), Real and Complex Analysis, 3rd Ed., McGraw-Hill.

[30] Soyster, A. L.(1973). Convex programming with set-inclusive constraints and applications to inexact linear programming. Operations Research, 21, 1154-1157.

[31] Talluri, K. T., van Ryzin, G. J. (2004). The Theory and Practice of Revenue Management. New York: Springer.

[32] Weatherford, L. R., Bodily, S. E. (1992). A taxonomy and research overview of perishable-asset revenue management: Yield management, overbooking and pricing. Operations Research, 40, 831-844.

[33] Zabel, E. (1970). Monopoly and uncertainty. The Review of Economic Studies, 37(2): 205-219.

[34] Zhang, Y. (2007). General Robust-Optimization Formulation for Nonlinear Programming. Journal of Optimization Theory and Applications, 132(1), 111-124.

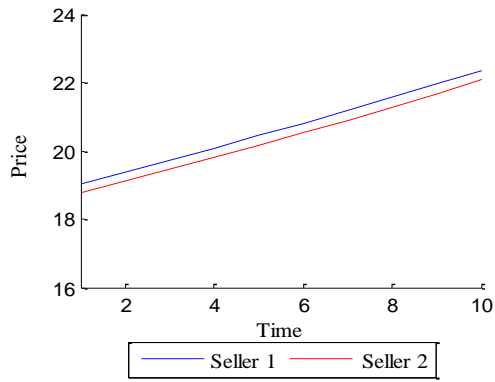

Figure 1: Price for each seller

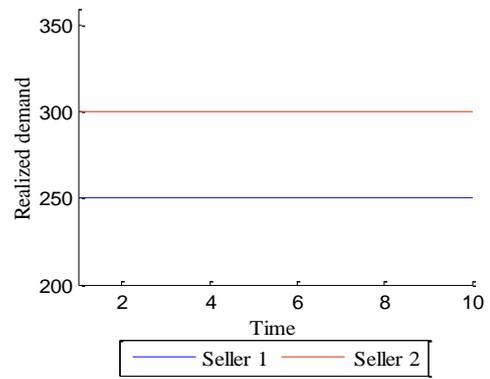

Figure 2: Realized demand for each seller

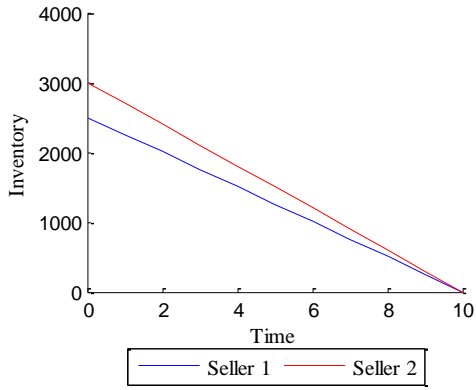
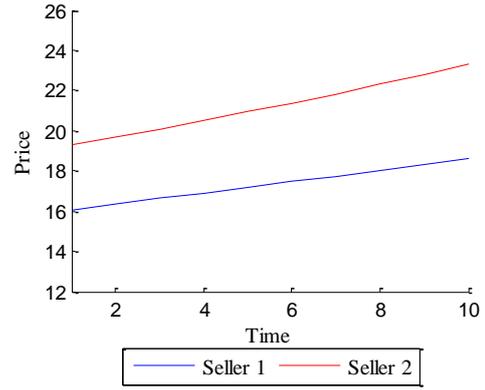

Figure 3: Inventory change for each seller    Figure 4: Price for each seller

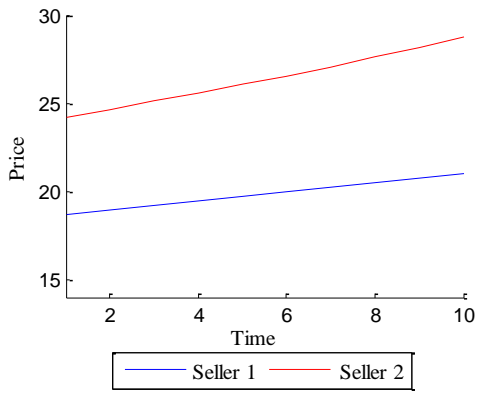
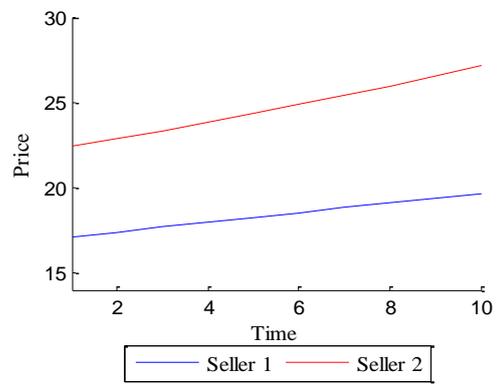

Figure 5: Price for each seller in Scenario 1    Figure 6: Price for each seller in Scenario 4

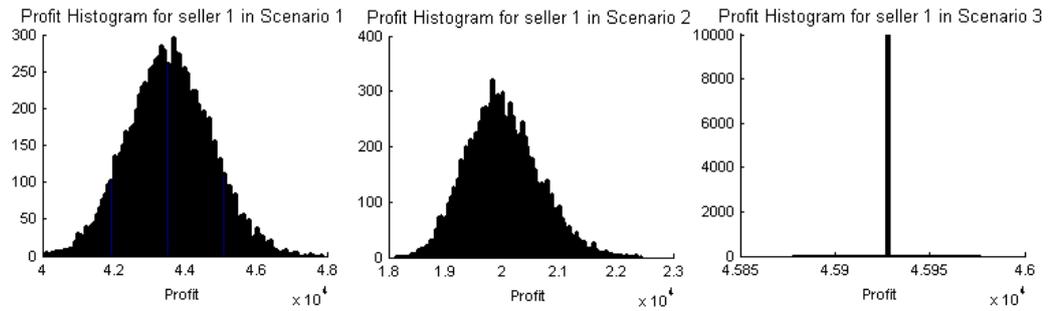

Figure 7: Profit histogram for seller 1 in Scenario 1-3 with beta (1, 3)

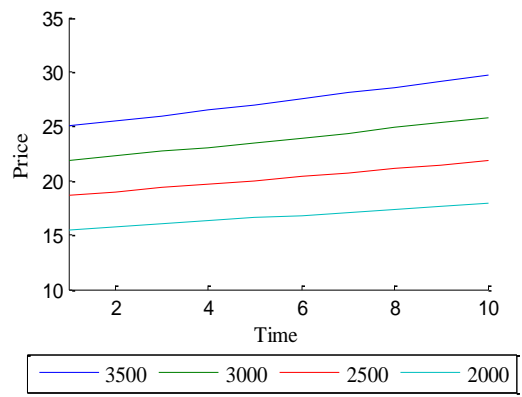
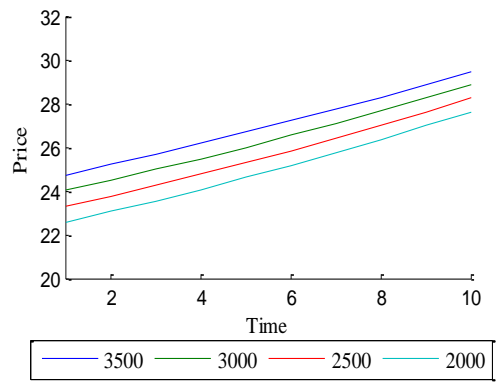

Figure 8: Prices for seller 1 with different value of $\alpha_1$

Figure 9: Prices for seller 2 with different value of $\alpha_1$

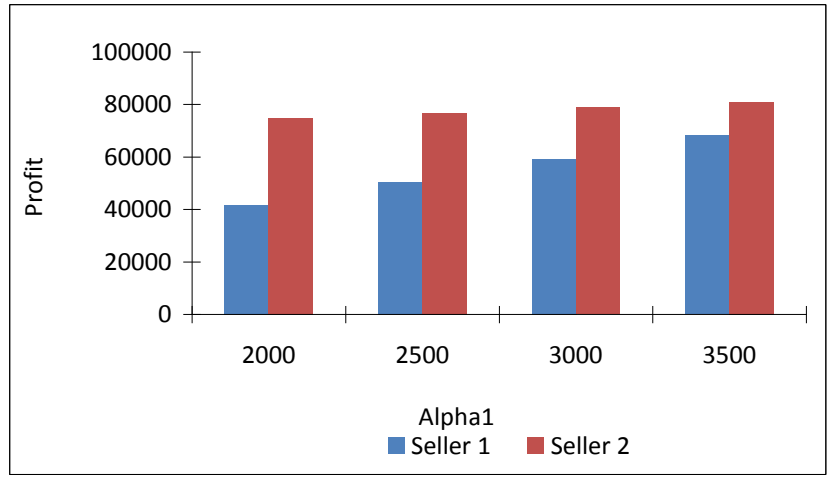

Figure 10: The effect of $\alpha_1$ on profits for each seller

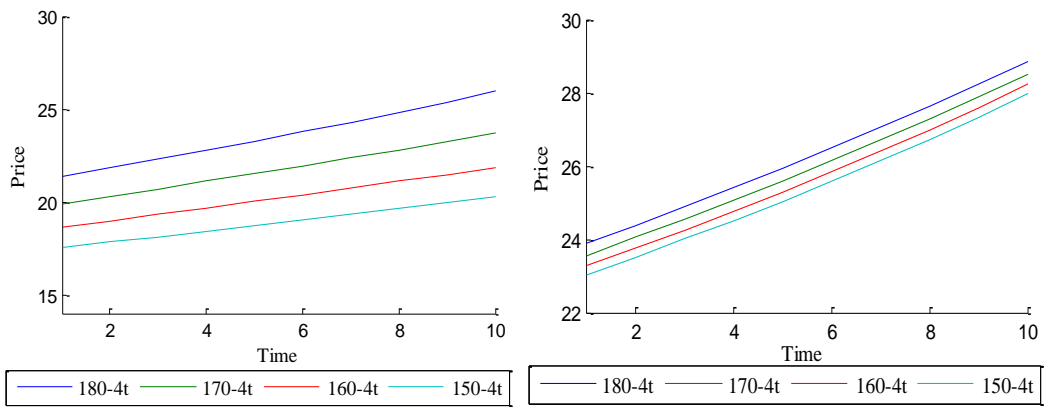

Figure 11: Prices for seller 1 with different value of $\beta_1$

Figure 12: Prices for seller 2 with different value of $\beta_1$

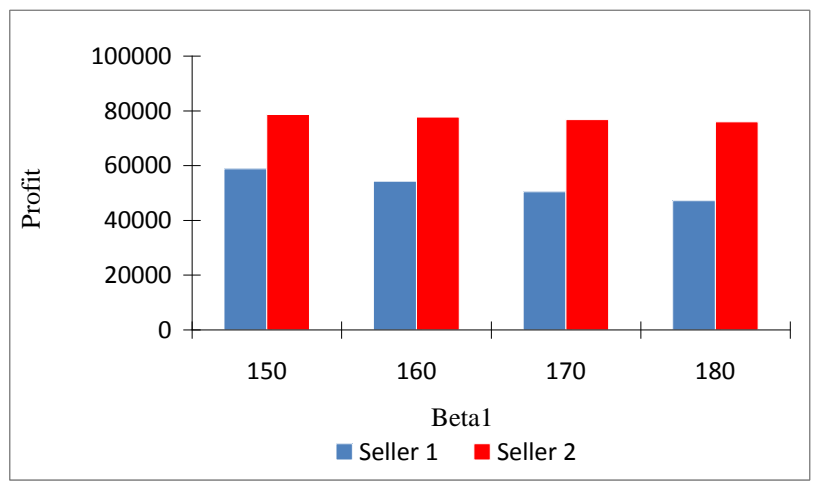

Figure 13: The effect of $\beta_1$ on profits for each seller

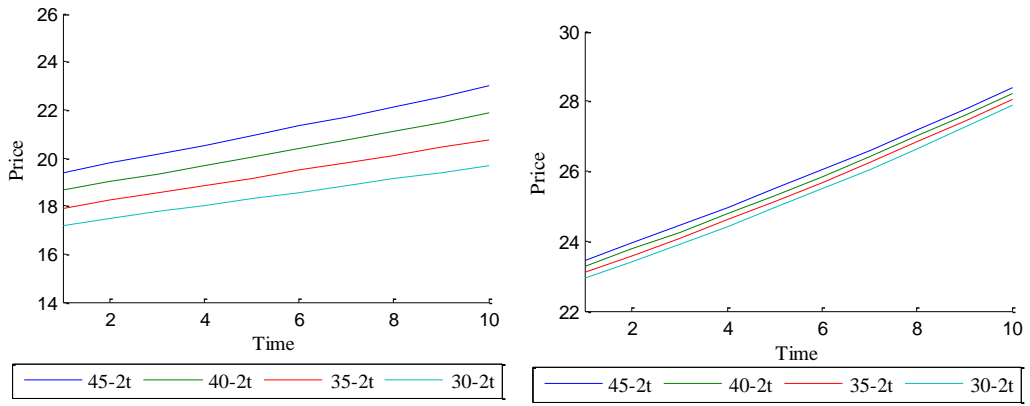

Figure 14: Prices for seller 1 with different value of $\gamma_2$

Figure 15: Prices for seller 2 with different value of $\gamma_2$

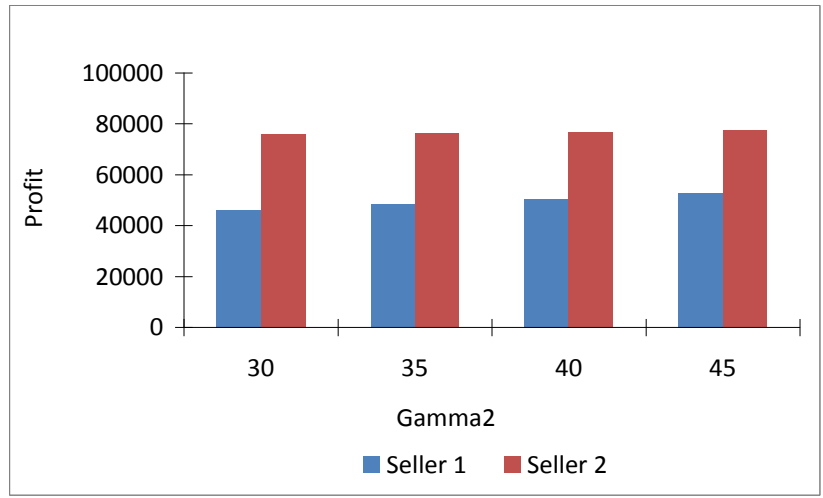

Figure 16: The effect of $\gamma_2$ on profits for each seller

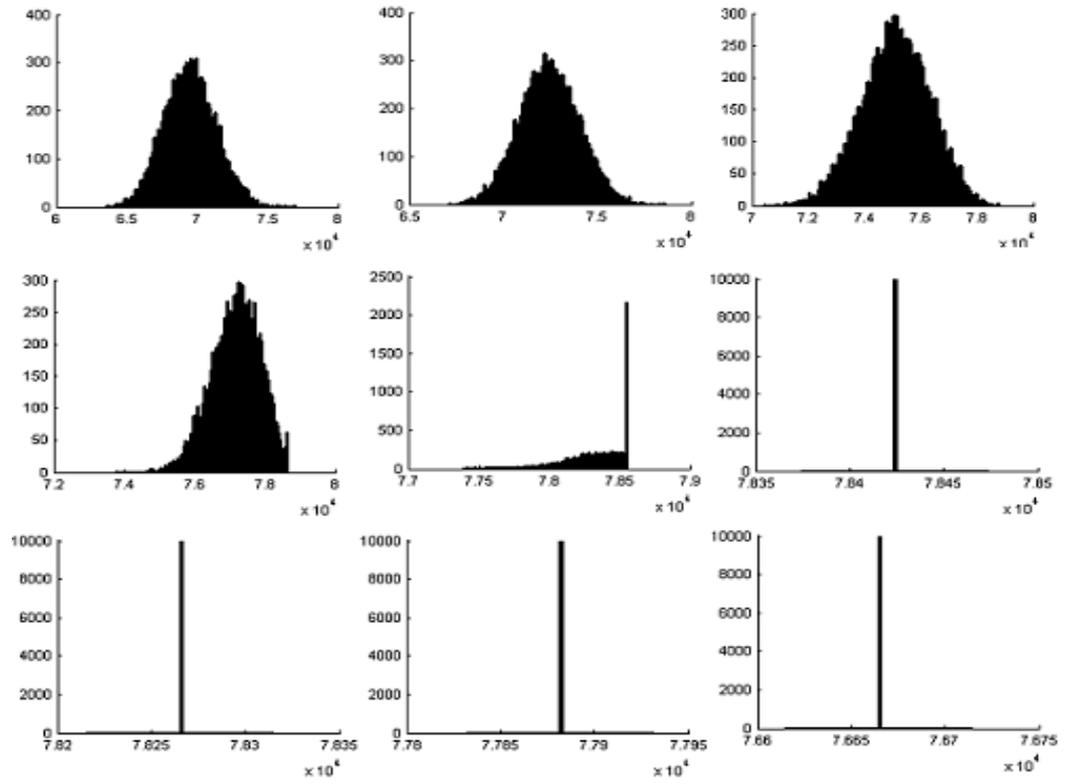

Figure 17: The profit histogram for seller 2 in Case I

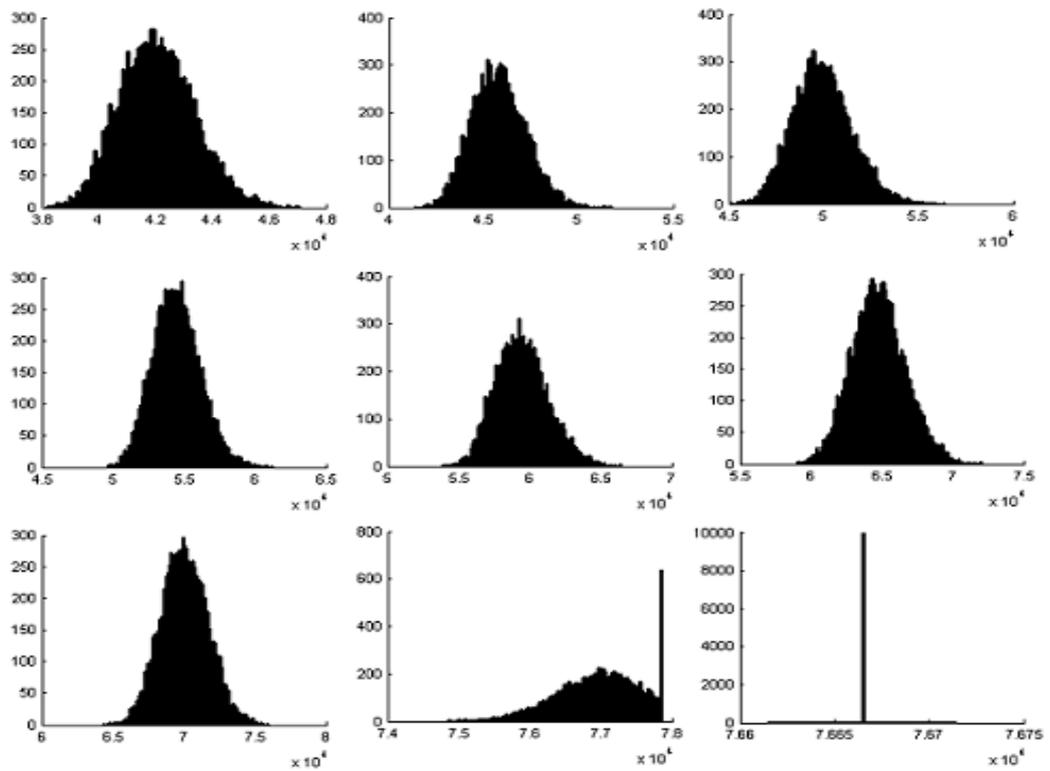

Figure 18: The profit histogram for seller 2 in Case II

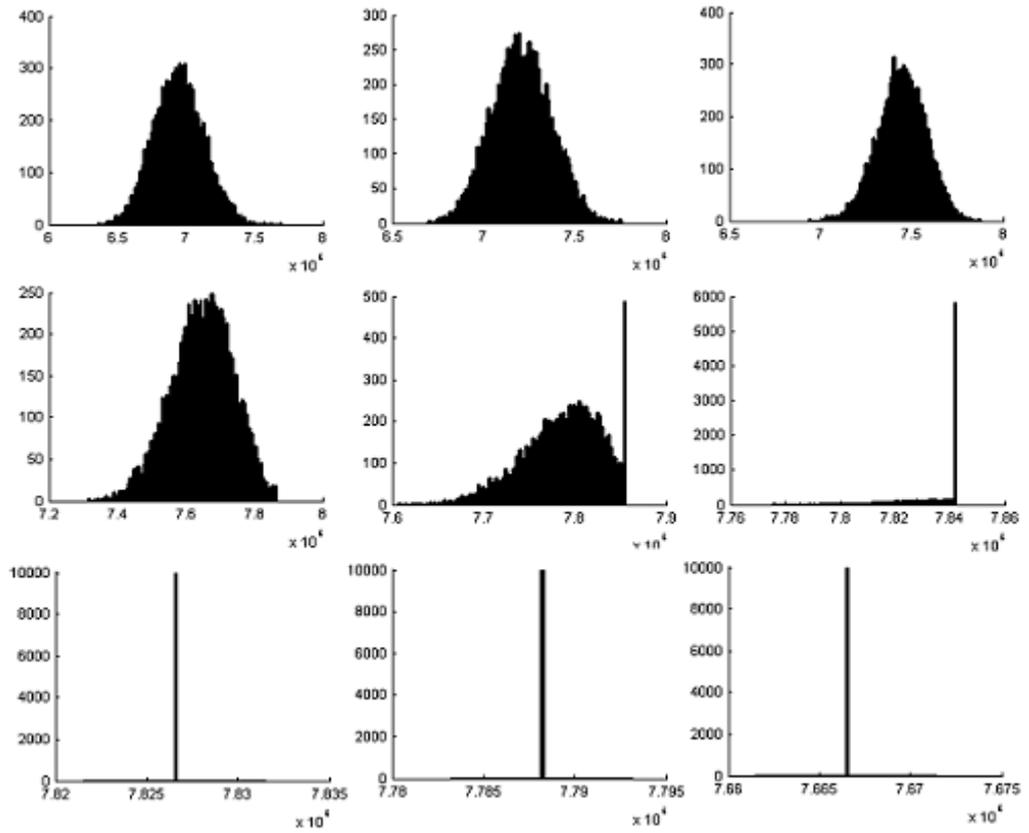

Figure 19: The profit histogram for seller 2 in Case III

Table 1 Profit for each seller

|  | Seller 1 | Seller 2 |
|---|---|---|
| Profit | 51,542 | 61,021 |

Table 2 Profit for each seller

|  | Seller 1 | Seller 2 |
|---|---|---|
| Profit | 43,238 | 63,542 |

Table 3 Pricing policies for both sellers under different scenarios

|  | Pricing policies | |
|---|---|---|
| Scenario | Seller 1 | Seller 2 |
| 1 | N | N |
| 2 | N | R |
| 3 | R | N |
| 4 | R | R |

Table 4 Range of profit for each seller with beta (1, 1) distribution

|  | Profit of seller 1 | | | | Profit of seller 2 | | | |
|---|---|---|---|---|---|---|---|---|
| Scenario | Min. | Max. | Av. | SD | Min. | Max. | Av. | SD |
| 1 | 42,271 | 49,533 | 46,730 | 1,145 | 66,247 | 79,003 | 74,528 | 1,794 |
| 2 | 18,962 | 25,882 | 22,510 | 932 | 73,885 | 73,885 | 73,885 | 0 |
| 3 | 45,928 | 45,928 | 45,928 | 0 | 40527 | 54,297 | 47,368 | 1,913 |
| 4 | 45,928 | 45,928 | 45,928 | 0 | 73,885 | 73,885 | 73,885 | 0 |

Table 5 Range of profit for each seller with beta (1, 3) distribution

|  | Profit of seller 1 | | | | Profit of seller 2 | | | |
|---|---|---|---|---|---|---|---|---|
| Scenario | Min. | Max. | Av. | SD | Min. | Max. | Av. | SD |
| 1 | 39,806 | 47,425 | 43,607 | 1,135 | 63,899 | 76,529 | 69,517 | 1,822 |
| 2 | 18,181 | 22,866 | 20,030 | 619 | 73,885 | 73,885 | 73,855 | 0 |
| 3 | 45,928 | 45,928 | 45,928 | 0 | 38,112 | 47,452 | 42,075 | 1,298 |
| 4 | 45,928 | 45,928 | 45,928 | 0 | 73,885 | 73,885 | 73,855 | 0 |